\documentclass{article}

\usepackage{amssymb}
\usepackage{amsthm}
\usepackage{amsmath}
\usepackage{bigints}
\usepackage{bm}
\usepackage{xcolor}

\usepackage{arxiv}

\usepackage[utf8]{inputenc} 
\usepackage[T1]{fontenc}    
\usepackage{hyperref}       
\usepackage{url}            
\usepackage{booktabs}       
\usepackage{amsfonts}       
\usepackage{nicefrac}       
\usepackage{microtype}      
\usepackage{lipsum}		
\usepackage{graphicx}
\usepackage{natbib}
\usepackage{doi}

\title{A general theoretical scheme for shape-programming of incompressible hyperelastic shells through differential growth
}

\author{
  Zhanfeng Li$^{1,3}$, Jiong Wang\thanks{Corresponding author. Email: \texttt{ctjwang@scut.edu.cn} Tel.: +86 13926459861, Fax: +86 21-87114460.} $^{,1,2}$, Mokarram Hossain$^{3}$, Chennakesava Kadapa$^{4}$ \\
  ${^1}$School of Civil Engineering and Transportation, South China University of Technology, Guangzhou, China\\
  ${^2}$State Key Laboratory of Subtropical Building Science, South China University of Technology, Guangzhou, China\\
  ${^3}$Zienkiewicz Centre for Computational Engineering (ZCCE), Swansea University, Swansea, United Kingdom\\
  ${^4}$School of Engineering and the Built Environment, Edinburgh Napier University, Edinburgh, United Kingdom\\
}

\renewcommand{\shorttitle}{A general theoretical scheme for shape-programming of shells}

\begin{document}
\maketitle

\begin{abstract}
	In this paper, we study the problem of shape-programming of incompressible hyperelastic shells through differential growth. The aim of the current work is to determine the growth tensor (or growth functions) that can produce the deformation of a shell to the desired shape. First, a consistent finite-strain shell theory is introduced. The shell equation system is established from the 3D governing system through a series expansion and truncation approach. Based on the shell theory, the problem of shape-programming is studied under the stress-free assumption. For a special case in which the parametric coordinate curves generate a net of curvature lines on the target surface, the sufficient condition to ensure the vanishing of the stress components is analyzed, from which the explicit expression of the growth tensor can be derived. In the general case, we conduct the variable changes and derive the total growth tensor by considering a two-step deformation of the shell. With these obtained results, a general theoretical scheme for shape-programming of thin hyperelastic shells through differential growth is proposed. To demonstrate the feasibility and efficiency of the proposed scheme, several nature-inspired examples are studied. The derived growth tensors in these examples have also been implemented in the numerical simulations to verify their correctness and accuracy. The simulation results show that the target shapes of the shell samples can be recovered completely. The scheme for shape-programming proposed in the current work is helpful in designing and manufacturing intelligent soft devices.
  \end{abstract}


\section{Introduction}
\label{sec:1}

Growth of soft biological tissues and swelling (or expansion) of soft polymeric gels are commonly observed in nature \citep{ambr2011,liu2015}. Due to the inhomogeneity or incompatibility of the growth fields, soft material samples usually exhibit diverse morphological changes and surface pattern evolutions during the growing processes, which is referred to as the `differential growth' and has attracted extensive research interest in recent years \citep{AlainGoriely2005,Li2012,kemp2014,Huang2018}. To fulfill the requirements of engineering applications, it is usually desired that the configurations of soft material samples are controllable during the growing processes, such that certain kinds of functions are realized. This goal can be achieved through sophisticated composition or architectural design in the soft material samples. The technique is known as `shape-programming' \citep{liu2016,vanManen2018}, which has been utilized for manufacturing a variety of intelligent soft devices, e.g., biomimetic 4D printing of flowers \citep{Gladman2016}, pressure-actuated deforming plate \citep{Siefert2019}, pasta with transient morphing effect \citep{Tao2021}, and polymorphic metal-elastomer composite \citep{Hwang2022}.

From the viewpoint of solid mechanics, soft materials can be treated as certain kinds of hyperelastic materials. The growth field in a soft material sample is usually modeled by incorporating a growth tensor. Due to the residual stresses triggered by the incompatibility of the growth field, as well as the external loads and boundary restrictions, the sample also undergoes elastic deformations. Thus, the total deformation gradient tensor should be decomposed into an elastic strain tensor and a growth tensor \citep{kond1987,Rodriguez1994,BenAmar2005}. The elastic incompressible constraint should also be adopted since the elastic deformations of soft materials are typically isochoric \citep{Wex2015, Kadapa2021}. Based on these constitutive and kinematic assumptions, the growth behaviors of soft material samples can be studied by solving the system of mechanical field equations. Because of the inherent nonlinearities in the large growth-induced deformations, mechanical instabilities can also be triggered in the soft material samples \citep{BenAmar2005,li2011,gori2017,Pezzulla2018,Xu2020}.

Despite the numerous studies on the growth behaviors of soft material samples, the majority of the modeling works pay attention to the direct problem, i.e., determining the deformations of soft material samples when the growth fields are specified. However, in order to utilize the shape-programming technique for engineering applications, one needs to study the inverse problem. That is, how to determine the growth fields in the samples such that the current configurations induced by differential growth can achieve any target shapes? This inverse problem has also been studied in some previous works \citep[cf.][]{dias2011,jone2015,acha2019,wang2019,nojo2021,Li2022,Wang2022}. In these works, the initial configurations of soft material samples usually have the thin plate form. Although the shell form is more common in nature and engineering fields, it is seldom chosen as the initial configuration of the soft material samples due to the difficulties associated with modelling shell structures.

To achieve the goal of shape-programming, a prerequisite is to predict the relations between the growth fields and the morphologies of soft material samples. It is thus of significance to establish an efficient and accurate mathematical model by taking configurations of samples, material properties, boundary conditions and other factors into account. In terms of shell theories for growth deformations, the Kirchhoff shell theory has been adopted to describe mechanical behavior in growing soft membranes \citep{Vetter2013,Rausch2014}, which relies on \textit{ad hoc} assumptions of the stress components and deformation gradient. Another shell theory is proposed based on the non-Euclidean geometry, where the deformation of samples is determined by the intrinsic geometric properties attached to surfaces, such as the first and second fundamental forms, and the applied growth fields \citep{Sadik2016,Pezzulla2018}. In \citet{Song2016}, a consistent finite-strain shell theory has been proposed within the framework of nonlinear elasticity, where the shell equation is derived from the 3D formulation through a series-expansion and truncation approach. To apply this theory for growth-induced deformations, \citet{Yu2022} incorporated the growth effect through the decomposition of the deformation gradient and derived the shell equation system for soft shell samples.

In this paper, we aim to propose a general theoretical scheme for shape-programming of incompressible hyperelastic shells through differential growth. Following the shell theory proposed in \citet{Yu2022}, the shell equation system is established from the 3D governing system, where a series expansion and truncation approach is adopted. To fulfill the purpose of shape-programming, the shell equation system is tackled by assuming that all the stress components vanish. Under this stress-free assumption, we first consider a special case in which the parametric coordinate curves generate a net of curvature lines on the target surface. By analyzing the sufficient condition to ensure the vanishing of the stress components, the explicit expression of the growth tensor is derived (i.e., the inverse problem is solved), which depends on the intrinsic geometric properties of the target surface. In the general case that the parametric coordinate curves cannot generate a net of curvature lines on the target surface, we conduct the variable changes and derive the total growth tensor by considering a two-step deformation of the shell sample. Based on these results, a theoretical scheme for shape-programming of hyperelastic shells is formulated. The feasibility and efficiency of this scheme are demonstrated by studying several typical examples.

This paper is organized as follows. In section 2, the finite-strain shell theory for modeling the growth behaviors of thin hyperelastic shells is introduced. In section 3, the problem of shape programming is solved and the theoretical scheme is proposed. In section 4, some typical examples are studied to show the efficiency of the theoretical scheme. Finally, some conclusions are drawn. In the following notations, the Greek letters $(\alpha, \beta, \gamma ...)$ run from 1 to 2, and the Latin letters $(i, j, k ...)$ run from 1 to 3. The repeated summation convention is employed and a comma preceding indices ${(\cdot)}_{,}$ represents the differentiation.


\section{The finite-strain shell theory}
\label{sec:2}

In this section, we first formulate the 3D governing system for modeling the growth behavior of a thin hyperelastic shell. Then, through a series-expansion and truncation approach, the finite-strain shell equation system of growth will be established.


\subsection{Kinematics and the 3D governing system}
\label{sec:2.1}

\begin{figure}[htbp]
    \centering
    \includegraphics[width=0.9\textwidth]{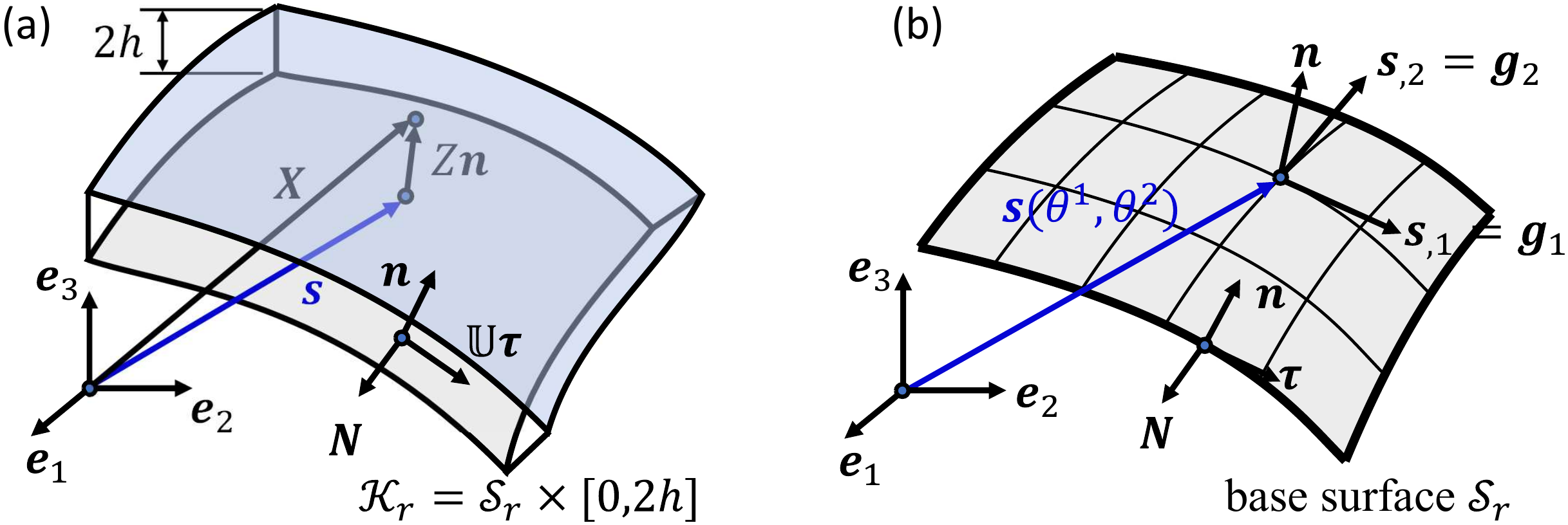}
    \caption{Position vector in the reference configuration $\mathcal{K}_r$: (a) reference configuration of the shell and decomposition of the position vector $\mathbf{X}$; (b) the curvilinear coordinate system and the local covariant basis on the base surface $\mathcal{S}_r$ of the shell.}
    \label{Fig:Configuration}
\end{figure}

We consider a thin homogeneous hyperelastic shell locating in the three-dimensional (3D) Euclidean space $R^3$. Within an orthonormal frame $\{O;\mathbf{e}_1,\mathbf{e}_2,\mathbf{e}_3\}$, the reference configuration of the shell occupies the region $\mathcal{K}_r=\mathcal{S}_r\times[0,2h]$, where the thickness parameter $h$ is much smaller than the dimensions of the base (bottom) surface $\mathcal{S}_r$ and its local radius of curvature. The position vector of a material point in the reference configuration $\mathcal{K}_r$ is denoted by $\mathbf{X}=X^i\mathbf{e}_i$ (cf. Fig. \ref{Fig:Configuration}(a)). The geometric description of a shell has been systematically reported in the literature \citep[cf.][]{Ciarlet2005,Steigmann2012,Song2016}, which is simply introduced below.

First, a curvilinear coordinate system $\{\theta^\alpha\}_{\alpha=1,2}$ is utilized to parametrize the base surface $\mathcal{S}_r$ of the shell in the reference configuration, which yields the parametric equation as
\begin{equation}
    \mathbf{s}(\theta^\alpha) = \left\{X^1(\theta^\alpha),X^2(\theta^\alpha),X^3(\theta^\alpha)\right\},\quad(\theta^\alpha)_{\alpha=1,2}\in\Omega_r.
    \label{Eq:ParaEqu}
\end{equation}
This parametric equation represents a continuous map from the region $\Omega_r \subset R^2$ to the surface $\mathcal{S}_r \subset R^3$. At a generic point on $\mathcal{S}_r$, the tangent vectors along the coordinate curves are given by $\mathbf{g}_\alpha=\mathbf{s}_{,\alpha}=\partial\mathbf{s}/\partial\theta^\alpha$, which span the tangent plane to the surface $\mathcal{S}_r$ at that point. The two vectors $\{\mathbf{g}_\alpha\}_{\alpha=1,2}$ are also referred to as the covariant basis of the tangent plane. Another two vectors $\{\mathbf{g}^\alpha\}_{\alpha=1,2}$ on the tangent plane can be determined unambiguously through the relations $\mathbf{g}_\alpha\cdot\mathbf{g}^\beta = \delta_{\alpha}^{\beta}$, which form the contravariant basis of the tangent plane. Then, the unit normal vector of the surface $\mathcal{S}_r$ should be defined by $\mathbf{n}=(\mathbf{g}_1\wedge\mathbf{g}_2) / \left|\mathbf{g}_1\wedge\mathbf{g}_2\right|$ (cf. Fig. \ref{Fig:Configuration}(b)). By denoting $\mathbf{g}_3=\mathbf{g}^3=\mathbf{n}$, $\{\mathbf{g}_i\}_{i=1,2,3}$ and $\{\mathbf{g}^i\}_{i=1,2,3}$ constitute two sets of right-handed orthogonal bases on the base surface $\mathcal{S}_r$. The first and second fundamental forms of the surface $\mathcal{S}_r$ can be written into \begin{equation}
    \mathrm{I}_r = g_{\alpha\beta}d\theta^\alpha d\theta^\beta,\ \ \ \mathrm{II}_r = b_{\alpha\beta}d\theta^\alpha d\theta^\beta,
    \label{Eq:1}
\end{equation}
where $g_{\alpha\beta}=\mathbf{g}_\alpha\cdot\mathbf{g}_\beta$ and $b_{\alpha\beta}=\mathbf{s}_{,\alpha\beta}\cdot\mathbf{n}$ are the fundamental quantities. Conventionally, the fundamental quantities are also denoted by
\begin{equation}
    \begin{aligned}
        &E_r = g_{11}, \quad F_r = g_{12} = g_{21}, \quad G_r = g_{22}, \\
        &L_r = b_{11}, \quad M_r = b_{12} = b_{21}, \quad N_r = b_{22}.
        \label{Eq:01}
    \end{aligned}
\end{equation}

As shown in Fig.\ref{Fig:Configuration}(a), the position vector $\mathbf{X}$ of a material point in the reference configuration $\mathcal{K}_r$ of the shell can be decomposed into
\begin{equation}
    \mathbf{X} = \mathbf{s}(\theta^\alpha) + Z \mathbf{n}(\theta^\alpha), \quad 0\leq Z\leq 2h,
    \label{Eq:2}
\end{equation}
where $Z$ is the coordinate of the point along the normal direction $\mathbf{n}$. Accordingly, the differential of $\mathbf{X}$ yields that
\begin{equation}
    d\mathbf{X}=d\mathbf{s}+Zd\mathbf{n}+\mathbf{n}dZ=\mathbf{g}_{\alpha} d\theta^{\alpha}+Z\mathbf{n}_{,\alpha}d\theta^{\alpha}+\mathbf{n}dZ.
    \label{Eq:3}
\end{equation}
From the Weingarten equation\citep{chenWH2017}, we have
\begin{equation}
    d\mathbf{n} = \mathbf{n}_{,\alpha}d\theta^{\alpha} = \left(\mathbf{n}_{,\alpha}\otimes\mathbf{g}^{\alpha}\right)\mathbf{g}_{\beta}d\theta^\beta = -\mathbb{K}d\mathbf{s},
    \label{Eq:4}
\end{equation}
where $\mathbb{K}=-\mathbf{n}_{,\alpha}\otimes\mathbf{g}^{\alpha}$ is the curvature tensor. The mean and Gaussian curvatures of the surface $\mathcal{S}_r$ are given by
\begin{equation}
    H = \frac{1}{2}\mathrm{tr}\left(\mathbb{K}\right), \quad K = \mathrm{Det}\left(\mathbb{K}\right).
    \label{Eq:5}
\end{equation}
By substituting (\ref{Eq:4}) into (\ref{Eq:3}), we obtain
\begin{equation}
    d\mathbf{X} = \mathbb{U}d\mathbf{s} + \mathbf{n}dZ = \hat{\mathbf{g}}_\alpha d\theta^\alpha + \mathbf{n}dZ,
    \label{Eq:6}
\end{equation}
where $\mathbb{U}=\mathbf{g}_{\alpha}\otimes\mathbf{g}^{\alpha}-Z\mathbb{K}$ and $\hat{\mathbf{g}}_\alpha=\mathbb{U}\mathbf{g}_\alpha$. We further denote $\hat{\mathbf{g}}^\alpha=\mathbb{U}^{-T}\mathbf{g}^\alpha$, then $\{\hat{\mathbf{g}}_\alpha\}_{\alpha=1,2}$ and $\{\hat{\mathbf{g}}^\alpha\}_{\alpha=1,2}$ form the covariant and contravariant base vectors at an arbitrary point in the shell, which are also orthogonal to $\mathbf{n}$. Notice that the thickness of the shell is much smaller than the radius of curvature of $\mathcal{S}_r$; thus, $\mathbb{U}$ should be an invertible tensor. From (\ref{Eq:6}), the area element on the base surface and the volume element in the shell can be written into
\begin{equation}
    \begin{aligned}
    & dA=\left|\mathbf{g}_1\wedge\mathbf{g}_2\right|d\theta^1d\theta^2=\sqrt{g_{11}g_{22}-g_{12}^2} \ d\theta^1d\theta^2,\\
    & dV=\mathrm{Det}(\mathbb{U})dAdZ=\left(1-2HZ+KZ^2\right)dAdZ.
    \end{aligned}
    \label{Eq:7}
\end{equation}
Regarding area element on the lateral surface $da$, the local differential follows \eqref{Eq:6} that
\begin{equation}
    \mathbf{N} da = (\mathbb{U}\bm{\tau}) \times \mathbf{n} \  ds dZ,
\end{equation}
where $\mathbf{N}$ is the outward normal unit vector of the lateral surface, and $\bm{\tau}$ is the unit tangent vector along the edge curve $\partial\mathcal{S}_r$ of the base surface, and $s$ is the arc-length variable on the edge curve $\partial\mathcal{S}_r$ of the base surface.
The norm of vector $(\mathbb{U}\bm{\tau})\times\mathbf{n}$ is denoted by $\sqrt{g_{\tau}}$ such that $ da = \sqrt{g_{\tau}} ds dZ$.

Due to the growth effect and the external loads, the configuration of the shell will deform from $\mathcal{K}_r$ to the current configuration $\mathcal{K}_t$ in $\mathcal{R}^3$. Within the orthonormal frame $\{O;\mathbf{e}_1,\mathbf{e}_2,\mathbf{e}_3\}$, the position vector of a material point in $\mathcal{K}_t$ is denoted by $\mathbf{x}(\theta^\alpha,Z)=x^i(\theta^\alpha,Z)\mathbf{e}_i$. The deformation gradient tensor $\mathbb{F}$ can then be calculated through
\begin{equation}
    \begin{aligned}
        \mathbb{F} = \mathbf{x}_{,\alpha} \otimes \hat{\mathbf{g}}^\alpha+\frac{\partial \mathbf{x}}{\partial Z} \otimes \mathbf{n} =\left(\nabla\mathbf{x}\right) \mathbb{U}^{-1}+\frac{\partial \mathbf{x}}{\partial Z} \otimes \mathbf{n},\\
    \end{aligned}
    \label{Eq:F}
\end{equation}
where $\nabla$ is the in-plane 2-D gradient on the base surface $\mathcal{S}_r$ ($\nabla\mathbf{x}=\mathbf{x}_{,\alpha} \otimes \mathbf{g}^\alpha$).

Following the basic assumption of growth mechanics \citep{kond1987,Rodriguez1994,BenAmar2005,GROH2022114839,DORTDIVANLIOGLU2017,MEHTA2021111026}, the deformation gradient tensor $\mathbb{F}$ is decomposed into
\begin{equation}
    \mathbb{F}=\mathbb{A}\mathbb{G},
    \label{Eq:FD}
\end{equation}
where $\mathbb{A}$ is the elastic strain tensor and $\mathbb{G}$ is the growth tensor. It is known that the rate of growth is relatively slow compared with the elastic response of the material, thus the distribution of the growth tensor $\mathbb{G}$ in the shell is assumed to be given and does not change.

As the elastic deformations of soft materials (e.g., soft biological tissues, polymeric gels) are generally isochoric, the following constraint equation should be adopted
\begin{equation}
    R(\mathbb{F},\mathbb{G})=J_GR_{0}(\mathbb{A})=J_G\left(\operatorname{Det}(\mathbb{A})-1\right)=0,
    \label{Eq:Incompressibility}
\end{equation}
where $J_G=\operatorname{Det}(\mathbb{G})$. Furthermore, we suppose the material has an elastic strain-energy function
\begin{equation}
    \phi (\mathbb{F},\mathbb{G}) = J_G \phi_0 (\mathbb{A})= J_G \phi_0 (\mathbb{FG}^{-1}).
    \label{Eq:Energyphi}
\end{equation}
Then, the nominal stress tensor $\mathbb{S}$ can be calculated through the constitutive equation
\begin{equation}
    \mathbb{S} = \frac{\partial \phi}{\partial \mathbb{F}} - p \frac{\partial R}{\partial \mathbb{F}} = J_G \mathbb{G}^{-1} \left( \frac{\partial \phi_0(\mathbb{A})}{\partial \mathbb{A}} -p \frac{\partial R_0 (\mathbb{A}) }{\partial \mathbb{A}}\right),
    \label{Eq:NominalStress1}
\end{equation}
where $p(\theta^\alpha,Z)$ is the Lagrange multiplier associated with the constraint \eqref{Eq:Incompressibility}.

During the growing process, the hyperelastic shell satisfies the following mechanical equilibrium equation
\begin{equation}
    \begin{aligned}
            \operatorname{Div} \mathbb{S} & = \left(\mathbb{S}_{,\alpha}\right)^T \hat{\mathbf{g}}^\alpha + \left(\frac{\partial \mathbb{S}}{\partial Z}\right)^T \mathbf{n} = \mathbf{0}, \quad \text {in}\ \ \mathcal{S}_r \times[0,2h]. \\
    \end{aligned}
    \label{Eq:DivS}
\end{equation}
We suppose that the bottom and top surfaces of the shell are subjected to the applied traction $\mathbf{q}^{\pm}$, which yields the boundary conditions
\begin{equation}
    \mathbb{S}^T\mathbf{n}|_{Z=0}=-\mathbf{q}^{-},\quad \mathbb{S}^T\mathbf{n}|_{Z=2h}=\mathbf{q}^{+}, \ \ \mathrm{on} \ \ \mathcal{S}_r.
    \label{Eq:Bou1}
\end{equation}
On the lateral surface $\partial\mathcal{S}_r\times[0,2h]$ of the shell, we suppose the applied traction is $\mathbf{q}(s,Z)$, where $s$ is the arc-length variable of boundary curve $\partial\mathcal{S}_r$. So, we also have the boundary condition
\begin{equation}
    \mathbb{S}^{T}\mathbf{N}=\mathbf{q}(s,Z) \quad \mathrm{on} \quad \partial\mathcal{S}_r\times[0,2h].
    \label{Eq:Bou2}
\end{equation}

Eqs. \eqref{Eq:Incompressibility} and \eqref{Eq:DivS} together with the boundary conditions \eqref{Eq:Bou1} and \eqref{Eq:Bou2} constitute the 3D governing system of the shell model, which contains the unknowns $\{ \mathbf{x},p\} $.


\subsection{Shell equation system}
\label{sec:2.2}

Starting from the 3D governing system of the shell model, the shell equation system can be derived through a series-expansion and a truncation approach. This approach has been proposed in \citet{Dai2014,Song2016,Wang.2016} for developing the consistent finite-strain plate and shell theories without the growth effect. In \citet{Wang.2018,Yu2022}, the finite-strain plate and shell theories of growth have also been established through this approach. For the sake of completeness of the current paper, the key steps of this approach to derive the shell equation system are introduced below (see \citet{Yu2022} for a comprehensive introduction). It should be noted that the derived shell equation system can attain the accuracy of $O(h^2)$. However, to fulfill the requirements of shape-programming in the following sections, we only need to present the shell equation to the asymptotic order of $O(h)$.

To eliminate the thickness variable $Z$ from the 3D governing system, we first conduct the series expansions of the unknowns as follows
\begin{equation}
    \begin{aligned}
            &\mathbf{x}(\theta^\alpha,Z) = \sum_{n=0}^{2} \frac{\mathbf{x}^{(n)}}{n!} Z^n + O \left(Z^3\right), \quad     p(\theta^\alpha,Z) = \sum_{n=0}^{2} \frac{p^{(n)}}{n!} Z^n + O \left(Z^3\right),
    \end{aligned}
    \label{Eq:Expansionxp}
\end{equation}
where $(\cdot)^{(n)}=\partial^{n}(\cdot) /\left.\partial Z^{n}\right|_{Z=0}$. According to \eqref{Eq:Expansionxp}, the deformation gradient tensor $\mathbb{F}$, the elastic strain tensor $\mathbb{A}$ and the nominal stress tensor $\mathbb{S}$ can also be expanded as
\begin{equation}
    \begin{aligned}
	\mathbb{F}&=\mathbb{F}^{(0)}+Z\mathbb{F}^{(1)}+O(Z^2), \\
	\mathbb{A}&=\mathbb{A}^{(0)}+Z\mathbb{A}^{(1)}+O(Z^2), \\
	\mathbb{S}&=\mathbb{S}^{(0)}+Z\mathbb{S}^{(1)}+O(Z^2).
    \end{aligned}
    \label{Eq:ExpansionFAS}
\end{equation}
Furthermore, we denote
\begin{equation}
    \begin{aligned}
    & \mathbb{G} = \mathbb{G}^{(0)} + Z \mathbb{G}^{(1)} + O(Z^2),\\
    & \mathbb{G}^{-1} = \bar{\mathbb{G}}^{(0)} + Z \bar{\mathbb{G}}^{(1)} + O(Z^2),\\
    & J_G\mathbb{G}^{-1} = \hat{\mathbb{G}}^{(0)} + Z \hat{\mathbb{G}}^{(1)} + O(Z^2).
    \end{aligned}
    \label{Eq:growthIn}
\end{equation}
Once the growth tensor $\mathbb{G}$ is given, $\bar{\mathbb{G}}^{(n)}$ and $\hat{\mathbb{G}}^{(n)}$ $(n=0,1)$ can be calculated directly.

By using the kinematic relations \eqref{Eq:F} and \eqref{Eq:FD}, the concrete expressions of $\mathbb{F}^{(n)}$ and $\mathbb{A}^{(n)}$ $(n=0,1)$ in terms of $\mathbf{x}^{(n)}$ $(n=0,1,2)$ can be derived. Further from the constitutive equation \eqref{Eq:NominalStress1}, we obtain
\begin{equation}
    \begin{aligned}
        \mathbb{S}^{(0)} &= \hat{\mathbb{G}}^{(0)} \left( \mathcal{A}^{(0)} - p^{(0)} \mathcal{R}^{(0)} \right), \\
        \mathbb{S}^{(1)} &= \hat{\mathbb{G}}^{(0)} \left( \mathcal{A}^{(1)} : \mathbb{A}^{(1)} -p^{(0)} \mathcal{R}^{(1)} : \mathbb{A}^{(1)} -p^{(1)} \mathcal{R}^{(0)} \right) + \hat{\mathbb{G}}^{(1)} \left( \mathcal{A}^{(0)} - p^{(0)} \mathcal{R}^{(0)} \right),
    \end{aligned}
    \label{Eq:S0S1JMPS}
\end{equation}
where $\mathcal{A}^{(n)}=\partial^{n+1}\phi_0/\partial\mathbb{A}^{n+1}|_{\mathbb{A} = \mathbb{A}^{(0)}}$ and $\mathcal{R}^{(n)}=\partial^{n+1} R_0/\partial\mathbb{A}^{n+1}|_{\mathbb{A} = \mathbb{A}^{(0)}}$ $(n=0,1)$.

We substitute \eqref{Eq:Expansionxp} and \eqref{Eq:ExpansionFAS} into the constraint equation \eqref{Eq:Incompressibility} and the mechanical equilibrium equation \eqref{Eq:DivS}. The coefficients of $Z^n$ $(n=0,1)$ in these equations should be zero, which yield that
\begin{equation}
    \begin{aligned}
       \mathrm{Det}\left( \mathbb{A}^{(0)} \right) -1 =0,\quad
       \mathcal{R}^{(0)} : \mathbb{A}^{(1)} = 0,
    \end{aligned}
    \label{Eq:ExpanDetA}
\end{equation}
and
\begin{equation}
    \begin{aligned}
        & \nabla \cdot \mathbb{S}^{(0)} + \left( \mathbb{S}^{(1)} \right)^T \mathbf{n}=\mathbf{0},\\
        & \nabla \cdot \mathbb{S}^{(1)} + \left( \mathbb{S}^{(2)} \right)^T \mathbf{n} + \mathbb{K}^T\mathbf{g}^\alpha\cdot\mathbb{S}^{(0)}_{,\alpha}=\mathbf{0}.
    \end{aligned}
    \label{Eq:ExpanDivS}
\end{equation}
Further substituting \eqref{Eq:Expansionxp} and \eqref{Eq:ExpansionFAS} into the boundary conditions \eqref{Eq:Bou1}, another two equations can be obtained
\begin{equation}
    \begin{aligned}
        &\left(\mathbb{S}^{(0)}\right)^T \mathbf{n} = -\mathbf{q}^-,\\
        &\left(\mathbb{S}^{(0)} + 2h \mathbb{S}^{(1)} + 2h^2\mathbb{S}^{(2)} \right)^T \mathbf{n} = \mathbf{q}^+.
    \end{aligned}
    \label{Eq:Bou3}
\end{equation}
Eqs. \eqref{Eq:ExpanDetA}$_2$ and \eqref{Eq:ExpanDivS}$_1$ constitute a linear system for $\mathbf{x}^{(2)}$ and $p^{(1)}$. By solving these two equations, we obtain \citep{Yu2022}
\begin{equation}
   \mathbf{x}^{(2)} = \mathbb{D}^{-1} \left( p^{(1)} \mathbf{y} -\mathbf{f} \right),\ \quad
   p^{(1)} = \frac{\mathbf{y} \cdot \mathbb{D}^{-1} \mathbf{f} - T}{\mathbf{y} \cdot \mathbb{D}^{-1} \mathbf{y}},
   \label{Eq:r2p1}
\end{equation}
where
\begin{equation*}
    \begin{aligned}
        (\mathbb{D})_{ij} =& \mathrm{Det}(\mathbb{G}^{(0)})\left( \mathcal{A}^{(1)} -p^{(0)} \mathcal{R}^{(1)} \right)_{kilj} \left( \left( \bar{\mathbb{G}}^{(0)} \right)^T \mathbf{n} \right)_k \left( \left( \bar{\mathbb{G}}^{(0)} \right)^T \mathbf{n} \right)_l,\\
        \mathbf{y} =& \mathcal{R}^{(0)^T} \hat{\mathbb{G}}^{(0)^T} \mathbf{n},\\
        \mathbf{f} =& \Big[ \left( \mathcal{A}^{(1)} - p^{(0)} \mathcal{R}^{(1)} \right) : \left[ \mathbb{F}^{(0)} \bar{\mathbb{G}}^{(1)} + \left( \mathbf{x}^{(1)} \otimes \nabla +\mathbf{x}^{(0)} \otimes \mathbb{K} \right) \bar{\mathbb{G}}^{(0)} \right] \Big]^T \hat{\mathbb{G}}^{(0)^T} \mathbf{n} \\
        &+ \left( \mathcal{A}^{(0)} - p^{(0)} \mathcal{R}^{(0)} \right)^T \hat{\mathbb{G}}^{(1)^T} \mathbf{n} + \nabla \cdot \mathbb{S}^{(0)}, \\
        T =& \mathrm{Det}(\mathbb{G}^{(0)}) \mathcal{R}^{(0)} : \left[ \left( \mathbf{x}^{(1)} \otimes \nabla + \mathbf{x}^{(0)} \otimes \nabla \mathbb{K} \right) \bar{\mathbb{G}}^{(0)} + \mathbb{F}^{(0)} \bar{\mathbb{G}}^{(1)} \right].
    \end{aligned}
\end{equation*}
The expressions of $\mathbf{x}^{(1)}$ and $p^{(0)}$ in terms of $\mathbf{x}^{(0)}$ can be obtained by solving the equations \eqref{Eq:ExpanDetA}$_1$ and \eqref{Eq:Bou3}$_1$. However, as these two equations are non-linear, the explicit expressions of $\mathbf{x}^{(1)}$ and $p^{(0)}$ can only be presented when a concrete form of the strain-energy function $\phi_0(\mathbb{A})$ is given.

To incorporate the effect of curvature of the shell, the factor $\mathrm{Det}(\mathbb{U})|_{Z=2h}=1-4hH+4h^2K$ is multiplied onto \eqref{Eq:Bou3}$_2$ \citep{Song2016}. In the remainder of this paper, we assume $H \leq O(1)$ and  $K \leq O(1)$, to ensure $1>|4hH|>|4h^2K|$ such that the terms consisting $h^2 H$ and $h^2 K$ can be dropped reasonably when the required order of equation is set as $O(h)$. By subtracting \eqref{Eq:Bou3}$_1$ from \eqref{Eq:Bou3}$_2$ and dividing it by $2h$, the following equation is obtained (where the terms of order higher than $O(h)$ have been dropped)
\begin{equation}
    \left(1-4h H \right) \left( \mathbb{S}^{(1)} \right)^T \mathbf{n} + h \left( \mathbb{S}^{(2)} \right)^T \mathbf{n} = \frac{\left(  1 - 4h H \right)  \mathbf{q}^+ + \mathbf{q}^-}{2h},\quad \mathrm{on}\ \mathcal{S}_r.
    \label{Eq:BcSubtractMu}
\end{equation}
By virtue of the relations given in \eqref{Eq:ExpanDivS}, \eqref{Eq:BcSubtractMu} can be rewritten into 2D vector shell equation
\begin{equation}
    \left(1-4h H\right) \nabla \cdot \mathbb{S}^{(0)} + h \left(\nabla \cdot \mathbb{S}^{(1)} + \mathbb{K}^T \mathbf{g}^\alpha \cdot \mathbb{S}^{(0)}_{,\alpha}\right) = -\frac{ \left(  1 - 4h H \right)  \mathbf{q}^+ + \mathbf{q}^- }{2h},\quad \mathrm{on}\ \mathcal{S}_r.
    \label{Eq:ShellEqs}
\end{equation}
which contains the unknown $\mathbf{x}^{(0)}(\theta^\alpha)$. In fact, $\mathbf{x}^{(0)}(\theta^\alpha)$ provides the parametric equation for the base surface $\mathcal{S}$ in the current configuration of the shell.

To establish a complete shell equation system, the boundary conditions on the edge $\partial\mathcal{S}_r$ should also be proposed. Based on the boundary condition \eqref{Eq:Bou2} in the 3D governing system, the following edge boundary conditions can be proposed
\begin{equation}
    \begin{aligned}
    & \left({\mathbb{S}^{(0)} + h \mathbb{S}^{(1)}}\right)^{T}\mathbf{N} = \int_{0}^{2h}\mathbf{q}(s,Z)dZ/(2h) = \bar{\mathbf{q}}, \\
    & \int_{\partial \mathcal{S}_r} \int_0^{2h} \left(\mathbb{S}^{T}\mathbf{N}\right) \wedge \left(\mathbf{x}(s,Z) - \mathbf{x}(s,h) \right) \sqrt{g_{\tau}} dZ ds \\
    & = \int_{\partial \mathcal{S}_r} \int_0^{2h} \mathbf{q}(s,Z) \wedge \left(\mathbf{x}(s,Z) - \mathbf{x}(s,h) \right) \sqrt{g_{\tau}} dZ ds = \bar{\mathbf{m}},
    \end{aligned}
    \label{Eq:ShellBcs}
\end{equation}
where $\bar{\mathbf{q}}$ and $\bar{\mathbf{m}}$ are the average traction and the bending moment (about the middle surface $Z=h$) applied on the lateral surface of the shell. Eqs. \eqref{Eq:ShellEqs} and \eqref{Eq:ShellBcs} constitute the shell equation system.


\section{Shape-programming of the thin hyperelastic shell}
\label{sec:3}

The shell equation system \eqref{Eq:ShellEqs}-\eqref{Eq:ShellBcs} can be applied to study the growth-induced deformations of the thin hyperelastic shell. For any given growth tensor $\mathbb{G}$ and boundary conditions, once the shell equation system is solved, the obtained solution $\mathbf{x}^{(0)}(\theta^\alpha)$ represents the base surface $\mathcal{S}$ of the shell in the current configuration $\mathcal{K}_t$.

The objective of the current work is to solve an inverse problem. That is, to ensure the shape of the base surface changes from $\mathcal{S}_r$ to a certain target shape $\mathcal{S}$, how to arrange the growth tensor (or growth functions) in the shell sample? This problem is referred to as `shape-programming' of thin hyperelastic shells \citep{liu2016}. For simplicity, we only consider the case that the surfaces of the shell are traction-free, i.e., $\mathbf{q}^\pm=\mathbf{q}(s,Z)=\mathbf{0}$ in \eqref{Eq:ShellEqs} and \eqref{Eq:ShellBcs}.

\begin{figure}[htbp]
    \centering
    \includegraphics[width=0.8\textwidth]{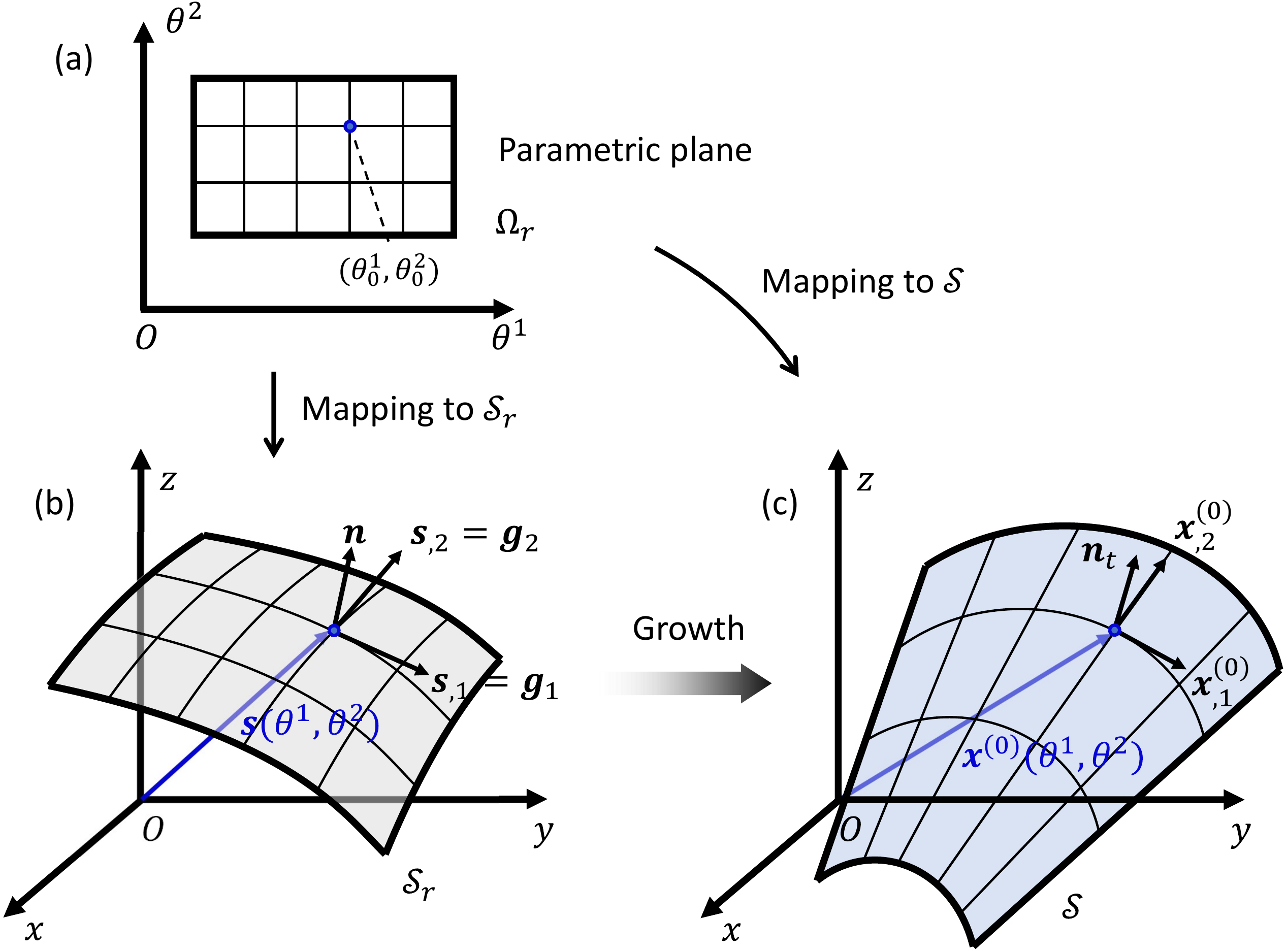}
    \caption{Illustration of the growth process and the mapping from parametric plane: (a) the variables region $\Omega_r$ on the parametric plane $\theta^1\theta^2$; (b) initial base surface $\mathcal{S}_r$; (c) target base surface $\mathcal{S}$.}
    \label{Fig:Kinematics}
\end{figure}

It should be noted that our goal is not to control the whole 3D configuration of the shell. As the shell equation system \eqref{Eq:ShellEqs}-\eqref{Eq:ShellBcs} is established on the base surface ($Z=0$) of the shell, we also focus on the base surface in solving the problem of shape-programming. The initial and current configurations of the base surface have the following parametric equations (as shown in Fig.\ref{Fig:Kinematics}):
\begin{equation}
    \begin{aligned}
        \mathcal{S}_r:&\ \mathbf{s}(\theta^\alpha) = \left\{X^1(\theta^\alpha),X^2(\theta^\alpha),X^3(\theta^\alpha) \right\},\quad (\theta^\alpha) \in \Omega_r,
    \end{aligned}
    \label{Eq:ParaSr}
\end{equation}
and
\begin{equation}
    \begin{aligned}
        \mathcal{S}:&\ \mathbf{x}^{(0)}(\theta^\alpha) = \left\{x^1(\theta^\alpha),x^2(\theta^\alpha),x^3(\theta^\alpha) \right\},\quad (\theta^\alpha) \in \Omega_r.
    \end{aligned}
    \label{Eq:ParaS}
\end{equation}
Eqs. \eqref{Eq:ParaSr} and \eqref{Eq:ParaS} can be viewed as two continuous mappings from the 2D region $\Omega_r$ to $\mathcal{S}_r$ and $\mathcal{S}$, respectively. By fixing one of the variables $\theta^1$ or $\theta^2$, variation of the the other variable can  generate the coordinate curves on the surfaces. All of these curves constitute the parametric curves net on $\mathcal{S}_r$ and $\mathcal{S}$. It has been introduced that on the initial configuration $\mathcal{S}_r$, the tangent vectors along the coordinate curves are $\mathbf{s}_{,\alpha}=\mathbf{g}_\alpha=\partial\mathbf{s}/\partial\theta^\alpha$ and unit normal vector is $\mathbf{n}=(\mathbf{g}_1\wedge\mathbf{g}_2) / \left|\mathbf{g}_1\wedge\mathbf{g}_2\right|$. Similarly, on base surface $\mathcal{S}$ in the current configuration, the tangent vectors along the coordinate curves are $\mathbf{x}^{(0)}_{,\alpha}=\partial\mathbf{x}^{(0)}/\partial\theta^\alpha$, and the unit normal vector is denoted by $\mathbf{n}_t= (\mathbf{x}_{,1}^{(0)}\wedge\mathbf{x}_{,2}^{(0)}) / \left|\mathbf{x}_{,1}^{(0)}\wedge\mathbf{x}_{,2}^{(0)}\right|$ (cf. Fig. \ref{Fig:Kinematics}). If $\mathcal{S}_r$ and $\mathcal{S}$ are regular surfaces, we always have $\mathbf{g}_1\wedge\mathbf{g}_2\neq0$ and $x_{,1}^{(0)}\wedge\mathbf{x}_{,2}^{(0)}\neq0$. Thus, the normal vector fields are well-defined on $\mathcal{S}_r$ and $\mathcal{S}$.

To facilitate the following derivations, we assume that the parametric curves net generated by $\{\mathbf{\theta}^\alpha\}$ is an orthogonal net of curvature lines on $\mathcal{S}_r$. This assumption means that the tangent vectors $\mathbf{g}_1$ and $\mathbf{g}_2$ are perpendicular to each other (i.e., $\mathbf{g}_1\cdot\mathbf{g}_2=0$) and they direct along the two principal directions at any point on $\mathcal{S}_r$. It is known that on a regular surface, such an orthogonal net always exists in the neighbour region of a non-umbilic point \citep{chenWH2017,topo2006}. Due to this assumption, some geometrical quantities defined in section \ref{sec:2.1} can be simplified into
\begin{equation}
    \begin{aligned}
        &F_r=g_{12}=g_{21}=0,\quad M_r=b_{12}=b_{21}=0,\\
        &\mathbb{K} = -\mathbf{n}_\alpha\otimes\mathbf{g}^\alpha = \kappa_1\mathbf{g}_1\otimes\mathbf{g}^1+\kappa_2\mathbf{g}_2\otimes\mathbf{g}^2, \\
        &H=\frac{1}{2}(\kappa_1+\kappa_2),\quad K=\kappa_1\kappa_2,\\
        &\mathbb{U} = \mathbb{I}_2-Z\mathbb{K} = (1-\kappa_1Z)\mathbf{g}_1\otimes\mathbf{g}^1+(1-\kappa_2Z)\mathbf{g}_2\otimes\mathbf{g}^2,
    \end{aligned}
    \label{Eq:GeoS}
\end{equation}
where $\kappa_1$ and $\kappa_2$ are called the principal curvatures.

With the above preparations, we begin to solve the problem of shape-programming of the thin hyperelastic shell. The major task is to reveal the relations between the growth tensor (or growth functions) and the geometrical properties of the target surface $\mathcal{S}$. Generally, the solution of shape-programming through differential growth may not be unique, i.e., the same target shape of the shell may be generated from different growth fields \citep{Wang2019c}. In this section, we focus on the case that the shell attains the stress-free state in the current configuration $\mathcal{K}_t$, i.e., all the components in $\mathbb{S}^{(0)}$ and $\mathbb{S}^{(1)}$ are zero. It is clear that in the stress-free condition, the shell equation system \eqref{Eq:ShellEqs} and \eqref{Eq:ShellBcs} is satisfied automatically.


\subsection{Growth tensor in a special case}
\label{Sec:3.1}

To analyze the relations between the growth tensor (or growth functions) and the geometrical properties of the target surface $\mathcal{S}$, we first assume that the coordinate curves of $\{\mathbf{\theta}^\alpha\}$ also generate a net of curvature lines on the base surface $\mathcal{S}$ in the current configuration. In this special case, the following specific form of the growth tensor will be adopted
\begin{equation}
    \begin{aligned}
    &\mathbb{G} = \mathbb{G}^{(0)} + Z\mathbb{G}^{(1)},\\
    &\mathbb{G}^{(0)} = \frac{\lambda_{1}^{(0)}}{\sqrt{E_r}} \mathbf{g}_1 \otimes \mathbf{g}^1 + \frac{\lambda_{2}^{(0)}}{\sqrt{G_r}}\mathbf{g}_2 \otimes \mathbf{g}^2+\mathbf{n} \otimes \mathbf{n},\\
    &\mathbb{G}^{(1)} = \frac{\lambda_{1}^{(1)}}{\sqrt{E_r}} \mathbf{g}_1 \otimes \mathbf{g}^1 + \frac{\lambda_{2}^{(1)}}{\sqrt{G_r}}\mathbf{g}_2 \otimes \mathbf{g}^2,
    \end{aligned}
    \label{Eq:growth}
\end{equation}
where $\lambda_{1}^{(0)}$, $\lambda_{2}^{(0)}$, $\lambda_{1}^{(1)}$ and $\lambda_{2}^{(1)}$ are the growth functions to be determined. By substituting \eqref{Eq:GeoS} and \eqref{Eq:growth} into the kinematic relations \eqref{Eq:F} and \eqref{Eq:FD}, the following expression of the elastic strain tensor $\mathbb{A}=\mathbb{F}\mathbb{G}^{-1}$ can be obtained
\begin{equation}
    \begin{aligned}
        \mathbb{A} = \frac{\sqrt{E_r}}{(1-\kappa_1Z)(\lambda_{1}^{(0)}+Z\lambda_{1}^{(1)})}\mathbf{x}_{,1}\otimes\mathbf{g}^1 +\frac{\sqrt{G_r}}{(1-\kappa_2Z)(\lambda_{2}^{(0)}+Z\lambda_{2}^{(1)})}\mathbf{x}_{,2}\otimes\mathbf{g}^2 + \frac{\partial \mathbf{x}}{\partial Z}\otimes\mathbf{n},
    \end{aligned}
    \label{Eq:FGS}
\end{equation}
The right Cauchy-Green strain tensor $\mathbb{C}=\mathbb{A}^T\mathbb{A}$ is then given by
\begin{equation}
    \begin{aligned}
        \mathbb{C} =& \frac{\mathbf{x}_{,1}\cdot\mathbf{x}_{,1}}{(1-\kappa_1Z)^2(\lambda_{1}^{(0)}+Z\lambda_{1}^{(1)})^2}{\hat{\mathbf{g}}}^1\otimes{\hat{\mathbf{g}}}^1 + \frac{\mathbf{x}_{,2}\cdot\mathbf{x}_{,2}}{(1-\kappa_2Z)^2(\lambda_{2}^{(0)}+Z\lambda_{2}^{(1)})^2}{\hat{\mathbf{g}}}^2\otimes{\hat{\mathbf{g}}}^2\\
        & + \frac{\mathbf{x}_{,1}\cdot\mathbf{x}_{,2}}{(1-\kappa_1Z)(1-\kappa_2Z)(\lambda_{1}^{(0)}+Z\lambda_{1}^{(1)})(\lambda_{2}^{(0)}+Z\lambda_{2}^{(1)})}\left({\hat{\mathbf{g}}}^1\otimes{\hat{\mathbf{g}}}^2 + {\hat{\mathbf{g}}}^2\otimes{\hat{\mathbf{g}}}^1\right)\\
        & + \frac{\mathbf{x}_{,1}\cdot\mathbf{x}_{,Z}}{(1-\kappa_1Z)(\lambda_{1}^{(0)}+Z\lambda_{1}^{(1)})}\left(\mathbf{n}\otimes{\hat{\mathbf{g}}}^1+{\hat{\mathbf{g}}}^1\otimes\mathbf{n}\right)\\
        & + \frac{\mathbf{x}_{,2}\cdot\mathbf{x}_{,Z}}{(1-\kappa_2Z)(\lambda_{2}^{(0)}+Z\lambda_{2}^{(1)})}\left(\mathbf{n}\otimes{\hat{\mathbf{g}}}^2+{\hat{\mathbf{g}}}^2\otimes\mathbf{n}\right) + \left(\mathbf{x}_{,Z}\cdot\mathbf{x}_{,Z}\right)\mathbf{n}\otimes\mathbf{n},
    \end{aligned}
    \label{Eq:CA}
\end{equation}
where ${\hat{\mathbf{g}}}^1=\sqrt{E_r}\mathbf{g}^1$ and ${\hat{\mathbf{g}}}^2=\sqrt{G_r}\mathbf{g}^2$ are two unit vectors. By substituting \eqref{Eq:Expansionxp}$_1$ into \eqref{Eq:CA} and conducting the series expansion of $\mathbb{C}$ with respect to $Z$, we have the following coefficients of $Z^0$ and $Z^1$
\begin{equation}
    \begin{aligned}
        \mathbb{C}^{(0)} =& \frac{\mathbf{x}^{(0)}_{,1}\cdot\mathbf{x}^{(0)}_{,1}}{{\lambda_{1}^{(0)}}^2}\hat{\mathbf{g}}^1\otimes\hat{\mathbf{g}}^1 + \frac{\mathbf{x}^{(0)}_{,2}\cdot\mathbf{x}^{(0)}_{,2}}{{\lambda_{2}^{(0)}}^2}\hat{\mathbf{g}}^2\otimes\hat{\mathbf{g}}^2 + \left(\mathbf{x}^{(1)}\cdot\mathbf{x}^{(1)}\right)\mathbf{n}\otimes\mathbf{n}\\
        &+ \frac{\mathbf{x}^{(0)}_{,1}\cdot\mathbf{x}^{(1)}}{\lambda_{1}^{(0)}} \left(\hat{\mathbf{g}}^1\otimes\mathbf{n} + \mathbf{n}\otimes\hat{\mathbf{g}}^1 \right) + \frac{\mathbf{x}^{(0)}_{,2}\cdot\mathbf{x}^{(1)}}{\lambda_{2}^{(0)}} \left(\hat{\mathbf{g}}^2\otimes\mathbf{n} + \mathbf{n}\otimes\hat{\mathbf{g}}^2 \right),
    \end{aligned}
    \label{Eq:CA0}
\end{equation}
\begin{equation}
    \begin{aligned}
        \mathbb{C}^{(1)} =& \frac{2}{{\lambda_{1}^{(0)}}^2}\left[\mathbf{x}^{(0)}_{,1}\cdot\mathbf{x}^{(1)}_{,1}+\left(\kappa_1-\frac{\lambda_{1}^{(1)}}{\lambda_{1}^{(0)}}\right)\mathbf{x}^{(0)}_{,1}\cdot\mathbf{x}^{(0)}_{,1}\right]\hat{\mathbf{g}}^1\otimes\hat{\mathbf{g}}^1\\
        & + \frac{2}{{\lambda_{2}^{(0)}}^2}\left[\mathbf{x}^{(0)}_{,2}\cdot\mathbf{x}^{(1)}_{,2}+\left(\kappa_2-\frac{\lambda_{2}^{(1)}}{\lambda_{2}^{(0)}}\right)\mathbf{x}^{(0)}_{,2}\cdot\mathbf{x}^{(0)}_{,2}\right]\hat{\mathbf{g}}^2\otimes\hat{\mathbf{g}}^2 + \left(2\mathbf{x}^{(1)}\cdot\mathbf{x}^{(2)}\right)\mathbf{n}\otimes\mathbf{n}\\
        & + \frac{\mathbf{x}^{(0)}_{,2}\cdot\mathbf{x}^{(1)}_{,1}+\mathbf{x}^{(0)}_{,1}\cdot\mathbf{x}^{(1)}_{,2}}{\lambda_{1}^{(0)}\lambda_{2}^{(0)}}\left(\hat{\mathbf{g}}^1\otimes\hat{\mathbf{g}}^2 + \hat{\mathbf{g}}^2\otimes\hat{\mathbf{g}}^1 \right) \\
        & + \left[\frac{\mathbf{x}^{(1)}_{,1}\cdot\mathbf{x}^{(1)}+\mathbf{x}^{(0)}_{,1}\cdot\mathbf{x}^{(2)}}{\lambda_{1}^{(0)}} + \frac{\mathbf{x}^{(0)}_{,1}\cdot\mathbf{x}^{(1)}}{\lambda_{1}^{(0)}}\left(\kappa_1-\frac{\lambda_{1}^{(1)}}{\lambda_{1}^{(0)}}\right)\right] \left(\hat{\mathbf{g}}^1\otimes\mathbf{n} + \mathbf{n}\otimes\hat{\mathbf{g}}^1 \right) \\
        & + \left[\frac{\mathbf{x}^{(1)}_{,2}\cdot\mathbf{x}^{(1)}+\mathbf{x}^{(0)}_{,2}\cdot\mathbf{x}^{(2)}}{\lambda_{2}^{(0)}} + \frac{\mathbf{x}^{(0)}_{,2}\cdot\mathbf{x}^{(1)}}{\lambda_{2}^{(0)}}\left(\kappa_2-\frac{\lambda_{2}^{(1)}}{\lambda_{2}^{(0)}}\right)\right] \left(\hat{\mathbf{g}}^2\otimes\mathbf{n} + \mathbf{n}\otimes\hat{\mathbf{g}}^2 \right),
    \end{aligned}
    \label{Eq:CA1}
\end{equation}
For isotropic incompressible hyperelastic material, it is known that the elastic strain-energy function only depends on the two invariants $I_1$ and $I_2$ of $\mathbb{C}$, i.e., $\phi_0(\mathbb{A})=\phi_0(I_1,I_2)$. Based on this constitutive form of $\phi_0$, the nominal stress tensor $\mathbb{S}$ given in \eqref{Eq:NominalStress1} can be rewritten into
\begin{equation}
    \begin{aligned}
        \mathbb{S} &= J_G \mathbb{G}^{-1} \left[ \frac{\partial \phi_0 }{\partial I_1} \frac{\partial I_1 }{\partial \mathbb{A}} + \frac{\partial \phi_0 }{\partial I_2} \frac{\partial I_2 }{\partial \mathbb{A}} -p \mathbb{A}^{-1} \right] \\
        &= J_G \mathbb{G}^{-1} \left[ 2 \frac{\partial \phi_0 }{\partial I_1} \mathbb{A}^T + 2 \frac{\partial \phi_0 }{\partial I_2} \left(I_1 \mathbb{I} - \mathbb{C}\right) \mathbb{A}^T -p \mathbb{A}^{-1}\right],
    \end{aligned}
    \label{Eq:NoS}
\end{equation}
where the relations $\partial I_1/\partial\mathbb{A}=2\mathbb{A}^T$ and $\partial I_2/\partial\mathbb{A} = 2(I_1\mathbb{I}-\mathbb{C})\mathbb{A}^T$ have been used. We denote
\begin{equation}
    \begin{aligned}
        & \mathbb{A}^{-1} = \bar{\mathbb{A}}^{(0)} + Z \bar{\mathbb{A}}^{(1)} + O(Z^2), \\
        & \mathbb{C} = \mathbb{C}^{(0)} + Z \mathbb{C}^{(1)} + O(Z^2), \quad I_1 = I_{1}^{(0)} + Z I_{1}^{(1)} + O(Z^2), \\
        & \frac{\partial \phi_0 }{\partial I_1} = d_{1}^{(0)} + Z d_{1}^{(1)} + O(Z^2), \quad \frac{\partial \phi_0 }{\partial I_2} = d_{2}^{(0)} + Z d_{2}^{(1)} + O(Z^2),\\
    \end{aligned}
\end{equation}
Then, the following expressions of $\mathbb{S}^{(0)}$ and $\mathbb{S}^{(1)}$ can be derived from \eqref{Eq:NoS}
\begin{equation}
    \begin{aligned}
        \mathbb{S}^{(0)} &= \hat{\mathbb{G}}^{(0)}\left\{ 2\left[d_{1}^{(0)}+d_{2}^{(0)}\left(I_{1}^{(0)} \mathbb{I}-\mathbb{C}^{(0)}\right)\right]\left(\mathbb{A}^{(0)}\right)^{T}-p^{(0)} \bar{\mathbb{A}}^{(0)}\right\} ,\\
        \mathbb{S}^{(1)} &= \hat{\mathbb{G}}^{(0)}\left\{2\left[d_{1}^{(0)}+d_{2}^{(0)}\left(I_{1}^{(0)} \mathbb{I}-\mathbb{C}^{(0)}\right)\right]\left(\mathbb{A}^{(1)}\right)^{T}\right\} \\
        &+\hat{\mathbb{G}}^{(0)}\left\{2\left[d_{1}^{(1)}+d_{2}^{(1)}\left(I_{1}^{(0)} \mathbb{I}-\mathbb{C}^{(0)}\right)+d_{2}^{(0)}\left(I_{1}^{(1)} \mathbb{I}-\mathbb{C}^{(1)}\right)\right]\left(\mathbb{A}^{(0)}\right)^{T}\right\} \\
        &+\hat{\mathbb{G}}^{(0)}\left(-p^{(0)} \bar{\mathbb{A}}^{(1)}-p^{(1)} \bar{\mathbb{A}}^{(0)}\right).
    \end{aligned}
    \label{Eq:S0S1}
\end{equation}
To ensure $\mathbb{S}^{(0)}$ and $\mathbb{S}^{(1)}$ to be zero tensors, one sufficient condition is that
\begin{equation}
    \begin{aligned}
        \mathbb{C}^{(0)} &= \mathbb{I}, \quad \mathbb{C}^{(1)} = \mathbf{0}, \quad p^{(0)} = 2\left(d_{1}^{(0)}+2 d_{2}^{(0)}\right), \quad p^{(1)} = 2\left(d_{1}^{(1)}+2 d_{2}^{(1)}\right).
    \end{aligned}
    \label{Eq:Sufficient}
\end{equation}
From \eqref{Eq:Sufficient}$_1$ and \eqref{Eq:Sufficient}$_2$, the growth functions $\{\lambda_{\alpha}^{(0)}\}$ and $\{\lambda_{\alpha}^{(1)}\}$ can be easily determined. In fact, by substituting \eqref{Eq:CA0} into \eqref{Eq:Sufficient}$_1$, we obtain
\begin{equation}
    \lambda_{1}^{(0)} =  \sqrt{E}, \quad \lambda_{2}^{(0)} = \sqrt{G},
    \label{Eq:GF1}
\end{equation}
where $E=\mathbf{x}^{(0)}_{,1}\cdot\mathbf{x}^{(0)}_{,1}$ and $G=\mathbf{x}^{(0)}_{,2}\cdot\mathbf{x}^{(0)}_{,2}$ are two of the first fundamental quantities of surface $\mathcal{S}$. As we assume the coordinate curves of $\{\mathbf{\theta}^\alpha\}$ generate a net of curvature lines on $\mathcal{S}$, another first fundamental quantity $F=\mathbf{x}^{(0)}_{,1}\cdot\mathbf{x}^{(0)}_{,2}=0$. Further from the relation \eqref{Eq:Sufficient}$_1$, we have
\begin{equation}
    \begin{aligned}
        &\mathbf{x}^{(0)}_{,1}\cdot\mathbf{x}^{(1)} = 0, \quad \mathbf{x}^{(0)}_{,2}\cdot\mathbf{x}^{(1)} = 0, \quad \mathbf{x}^{(1)}\cdot\mathbf{x}^{(1)} = 1, \\
        \Rightarrow &\  \mathbf{x}^{(1)} = \frac{\mathbf{x}^{(0)}_{,1}\wedge\mathbf{x}^{(0)}_{,2}}{\left|\mathbf{x}^{(0)}_{,1}\wedge\mathbf{x}^{(0)}_{,2}\right|} = \mathbf{n}_t.
    \end{aligned}
    \label{Eq:GF2}
\end{equation}
By substituting \eqref{Eq:GF1} and \eqref{Eq:GF2} into \eqref{Eq:CA1}, then from \eqref{Eq:Sufficient}$_2$, we obtain
\begin{equation}
    \lambda_1^{(1)} = \left(\kappa_1 - \frac{L}{E}\right)\sqrt{E},\quad \lambda_2^{(1)} = \left(\kappa_2 - \frac{N}{G}\right)\sqrt{G},
    \label{Eq:GF3}
\end{equation}
where $L=-\mathbf{x}^{(0)}_{,1}\cdot\mathbf{n}_{t,1}$ and $N=-\mathbf{x}^{(0)}_{,2}\cdot\mathbf{n}_{t,2}$ are two of the second fundamental quantities of surface $\mathcal{S}$. Another second fundamental quantity $M=-\mathbf{x}^{(0)}_{,1}\cdot\mathbf{n}_{t,2}=-\mathbf{x}^{(0)}_{,2}\cdot\mathbf{n}_{t,1}=0$ due to the net of curvature lines on $\mathcal{S}$. To ensure all the components of $\mathbb{C}^{(1)}$ to be zero, we also need to set $\mathbf{x}^{(2)}=\mathbf{0}$. By substituting \eqref{Eq:GF1} and \eqref{Eq:GF3} into \eqref{Eq:growth}, we obtain
\begin{equation}
    \begin{aligned}
        \mathbb{G} = & \left[1+Z \left(\kappa_1-\frac{L}{E}\right)\right]\sqrt{\frac{E}{E_r}} \mathbf{g}_1\otimes\mathbf{g}^1 \\
        &+ \left[1+Z \left(\kappa_2-\frac{N}{G}\right)\right]\sqrt{\frac{G}{G_r}} \mathbf{g}_2\otimes\mathbf{g}^2 + \mathbf{n} \otimes \mathbf{n},
    \end{aligned}
    \label{Eq:GT1}
\end{equation}
which is just the growth tensor that can result in the shape change of the base surface of the shell from $\mathcal{S}_r$ to $\mathcal{S}_t$ in the special case (i.e., the coordinate curves of $\{\mathbf{\theta}^\alpha\}$ generate a net of curvature lines on $\mathcal{S}$).

The growth functions $\lambda_{1}^{(0)}$ and $\lambda_{2}^{(0)}$ given in \eqref{Eq:GF1} have the same expressions as those obtained from the plate model (where the incompressible Neo-Hookean material is taken into account) \citep{Wang2022}. In fact, $\lambda_{1}^{(0)}$ and $\lambda_{2}^{(0)}$ just represent the extension or shrinkage of the material along the coordinate curves of $\{\theta^\alpha\}$ on $\mathcal{S}_r$. The growth functions $\lambda_{1}^{(1)}$ and $\lambda_{2}^{(1)}$ given in \eqref{Eq:GF3} involve the principal curvatures $\kappa_1$ and $\kappa_2$ of $\mathcal{S}_r$, which are different from the results of the plate model \citep{Wang2022}. It should be noted that the growth functions given in \eqref{Eq:GF1} and \eqref{Eq:GF3} are independent of the strain-energy function $\phi_0$, which should be valid for different kinds of hyperelastic shells. If the shell is made of incompressible Neo-Hookean material, the results \eqref{Eq:GF1} and \eqref{Eq:GF3} can be derived through another approach, which is introduced in \ref{app:1}.


\subsection{Growth tensor in general cases}
\label{Sec:3.2}

The formulas \eqref{Eq:GF1} and \eqref{Eq:GF3} are derived based on the assumption that the coordinate curves of variables $\{\theta^{\alpha}\}$ constitute an orthogonal net of curvature lines in the current configuration of the base surface $\mathcal{S}$. Generally, this assumption cannot be satisfied by the parametric equation $\mathbf{x}^{(0)}(\theta^\alpha)$. To tackle the problem in general cases, some further manipulations are required.

\begin{figure}[htbp]
    \centering
    \includegraphics[width=0.99\textwidth]{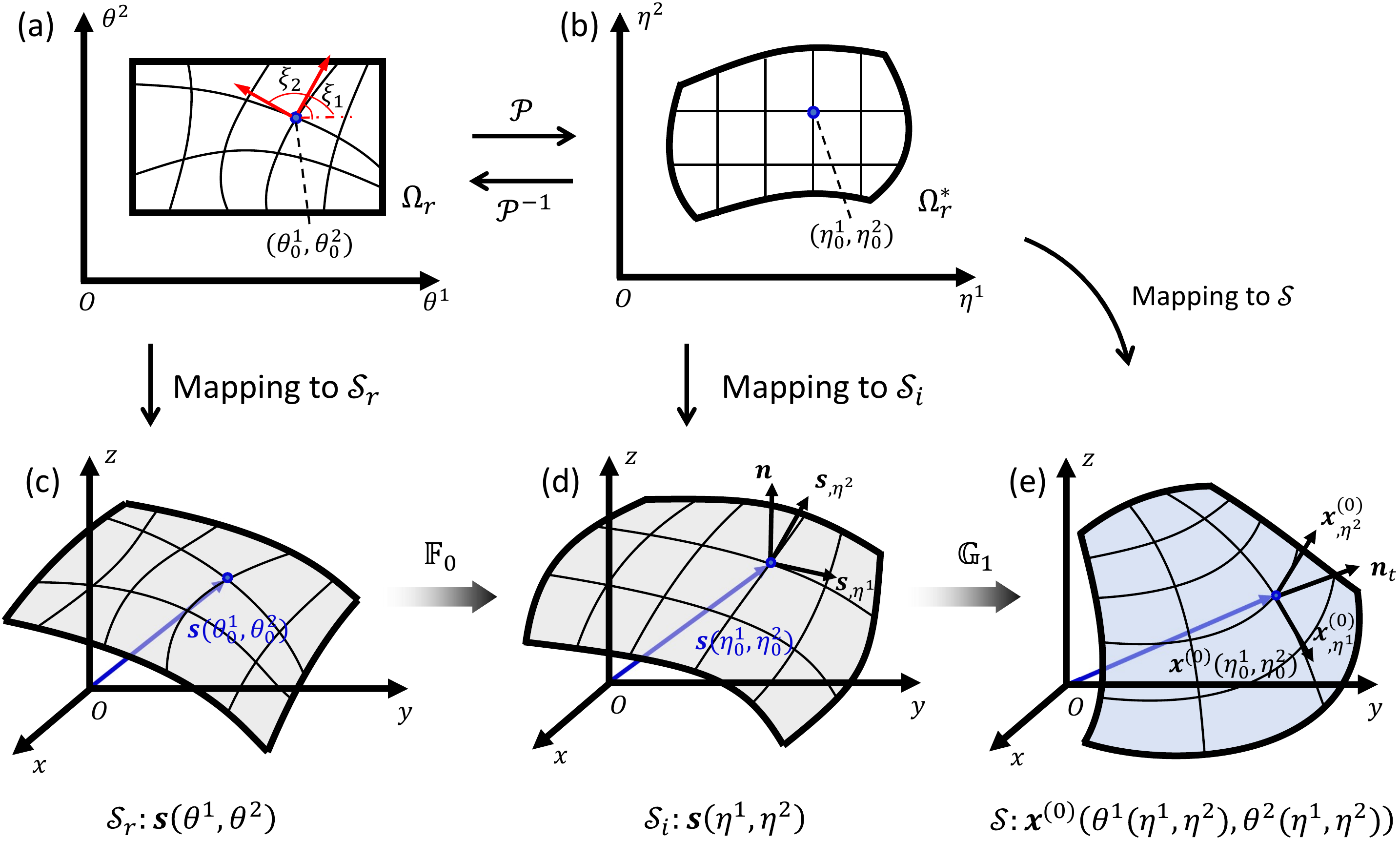}
    \caption{Bijection between the region $\Omega_r$ and $\Omega_r^*$, and the decomposition of the growth process: (a) the original region $\Omega_r$ on the parametric plane $\theta^1 \theta^2$; (b) the new region $\Omega_r^*$ on the parametric plane $\eta^1 \eta^2$; (c) the base surface $\mathcal{S}_r$ in the referential configuration $\mathcal{K}_r$; (d) the base surface $\mathcal{S}_i$ on the intermediate configuration $\mathcal{K}_i$; (e) the target base surface $\mathcal{S}$ in the current configuration $\mathcal{K}_t$.}
    \label{Fig:ChangeVariables}
\end{figure}

First, to generate a net of curvature lines on the surface $\mathcal{S}$, we consider the following change of variables
\begin{equation}
    \theta^1=\theta^1\left( \eta^1, \eta^2 \right), \quad \theta^2=\theta^2\left( \eta^1, \eta^2 \right),
    \label{Eq:Transformation1}
\end{equation}
where $\theta^1\left( \eta^1, \eta^2 \right)$ and $\theta^2\left( \eta^1, \eta^2 \right)$ are supposed to be sufficient smooth and the Jacobi determinant $\partial(\theta^1,\theta^2)/\partial(\eta^1,\eta^2)>0$. The transformation \eqref{Eq:Transformation1} defines a bijection between the planar parametric region $\Omega_r$ in the $\theta^1\theta^2$-plane and the planar parametric region $\Omega_r^*$ in the $\eta^1\eta^2$-plane (cf. Fig. \ref{Fig:ChangeVariables}).

Through the change of variables, the surface $\mathcal{S}$ has a new parametric equation $\mathbf{x}^{*}(\eta^1,\eta^2)=\mathbf{x}^{(0)} \left( \theta^1\left( \eta^1, \eta^2 \right),\theta^2\left( \eta^1, \eta^2 \right) \right)$. The first-order derivatives of $\mathbf{x}^{*}(\eta^1,\eta^2)$ are given by
\begin{equation}
    \begin{aligned}
        &\mathbf{x}_{, \eta^{1}}^{*} = \mathbf{x}_{, \theta^{1}}^{*} \frac{\partial \theta^{1}}{\partial \eta^{1}}+\mathbf{x}_{, \theta^{2}}^{*} \frac{\partial \theta^{2}}{\partial \eta^{1}}=A_{1}\left(\mathbf{x}_{, \theta^{1}}^{*} \cos \xi_{1}+\mathbf{x}_{, \theta^{2}}^{*} \sin \xi_{1}\right), \\
        &\mathbf{x}_{, \eta^{2}}^{*} = \mathbf{x}_{, \theta^{1}}^{*} \frac{\partial \theta^{1}}{\partial \eta^{2}}+\mathbf{x}_{, \theta^{2}}^{*} \frac{\partial \theta^{2}}{\partial \eta^{2}}=A_{2}\left(\mathbf{x}_{, \theta^{1}}^{*} \cos \xi_{2}+\mathbf{x}_{, \theta^{2}}^{*} \sin \xi_{2}\right),
        \end{aligned}
\end{equation}
where
\begin{equation}
    \begin{aligned}
        A_{1} & =\sqrt{\left(\frac{\partial \theta^{1}}{\partial \eta^{1}}\right)^{2}+\left(\frac{\partial \theta^{2}}{\partial \eta^{1}}\right)^{2}},\quad \cos \xi_{1}=\frac{1}{A_{1}}\frac{\partial \theta^{1}}{\partial \eta^{1}}, \quad \sin \xi_{1}= \frac{1}{A_{1}}\frac{\partial \theta^{2}}{\partial \eta^{1}}, \\
        A_{2} & =\sqrt{\left(\frac{\partial \theta^{1}}{\partial \eta^{2}}\right)^{2}+\left(\frac{\partial \theta^{2}}{\partial \eta^{2}}\right)^{2}},\quad \cos \xi_{2}= \frac{1}{A_{2}}\frac{\partial \theta^{1}}{\partial \eta^{2}}, \quad \sin \xi_{2}=\frac{1}{A_{2}}\frac{\partial \theta^{2}}{\partial \eta^{2}}.
        \end{aligned}
        \label{Eq:A1A2xi1xi2}
\end{equation}
To ensure that the new coordinate curves (i.e., the $\eta^1$- and $\eta^2$-curves on $\mathcal{S}$) constitute an orthogonal net of curvature lines, $\mathrm{x}_{, \eta^1}^{*}$ and $\mathrm{x}_{, \eta^2}^{*}$ should be aligned with the principal directions at any point on $\mathcal{S}$, which requires that the following equation is satisfied \citep{chenWH2017,topo2006}
\begin{equation}
    (L F-M E) \cos ^{2} \xi+(L G-N E) \cos \xi \sin \xi+(M G-N F) \sin ^{2} \xi = 0,
\end{equation}
where $\{E, F, G\}$ and $\{L, M, N\}$ are the first and second fundamental quantities of the surface $\mathcal{S}$ calculated from the original
parametric equation $\mathbf{x}^{(0)}(\theta^\alpha)$. On the other hand, as the Jacobi determinant $\partial(\theta^1,\theta^2)/\partial(\eta^1,\eta^2)>0$, we have the inverse Jacobi matrix
\begin{equation}
    \left(\begin{array}{ll}
    \frac{\partial \eta^{1}}{\partial \theta^{1}} & \frac{\partial \eta^{2}}{\partial \theta^{1}} \\
    \frac{\partial \eta^{1}}{\partial \theta^{2}} & \frac{\partial \eta^{2}}{\partial \theta^{2}}
    \end{array}\right)=\left(\begin{array}{ll}
    \frac{\partial \theta^{1}}{\partial \eta^{1}} & \frac{\partial \theta^{1}}{\partial \eta^{2}} \\
    \frac{\partial \theta^{2}}{\partial \eta^{1}} & \frac{\partial \theta^{2}}{\partial \eta^{2}}
    \end{array}\right)^{-1}.
    \label{Eq:JacobiMatrix}
\end{equation}
By virtue of \eqref{Eq:JacobiMatrix}, the differential forms $d\eta^{1}$ and $d\eta^{2}$ can be written into
\begin{equation}
    \begin{aligned}
    &d \eta^{1}=\frac{\partial \eta^{1}}{\partial \theta^{1}} d \theta^{1}+\frac{\partial \eta^{1}}{\partial \theta^{2}} d \theta^{2}=A_{1}^{*}\left(\sin \xi_{2} {d} \theta^{1}-\cos \xi_{2} {d} \theta^{2}\right), \\
    &d \eta^{2}=\frac{\partial \eta^{2}}{\partial \theta^{1}} d \theta^{1}+\frac{\partial \eta^{2}}{\partial \theta^{2}} d \theta^{2}=A_{2}^{*}\left(\cos \xi_{1} {d} \theta^{2}-\sin \xi_{1} {d} \theta^{1}\right),
    \end{aligned}
    \label{Eq:eta1eta2}
\end{equation}
where
\begin{equation}
    A_{1}^{*}=\frac{1}{{A}_{1} \left(\cos \xi_{1} \sin \xi_{2}- \sin \xi_{1} \cos \xi_{2}\right)},\quad A_{2}^{*}=\frac{1}{{A}_{2}\left( \cos \xi_{1} \sin \xi_{2}- \sin \xi_{1} \cos \xi_{2}\right)}.
\end{equation}
To make the differential forms $d \eta^1$ and $d \eta^2$ given in \eqref{Eq:eta1eta2} to be integrable, it necessary to derive the explicit expressions of the transformation between $\{\theta^1,\theta^2\}$ and $\{\eta^1,\eta^2\}$. To our knowledge, there are still no universal formulas that can be used to determine the integrating factors for any differential forms \citep{chenWH2017}. In some particular situations, the integrating factors can be obtained by adopting appropriate methods. Once the integrating factors are found, the explicit expressions of $\eta^{1}(\theta^1,\theta^2)$ and $\eta^{2}(\theta^1,\theta^2)$ can be obtained by the first integrals of the differential forms \eqref{Eq:eta1eta2}. Accordingly, the expressions of $\theta^{1}(\eta^1,\eta^2)$ and $\theta^{2}(\eta^1,\eta^2)$ are also obtained.

On the parametric variable region $\Omega_r^*$ in the $\eta^1\eta^2$-plane, we define a new surface $\mathcal{S}_i$ in $\mathcal{R}^3$, which has the following parametric equation
\begin{equation}
\mathcal{S}_i:\ \mathbf{s}(\eta^\alpha) = \left\{X^1(\eta^\alpha),X^2(\eta^\alpha),X^3(\eta^\alpha) \right\},\quad (\eta^\alpha) \in \Omega_r^*.
\label{Eq:ParaESi1}
\end{equation}
Notice that $\mathcal{S}_i$ and $\mathcal{S}_r$ have the same parametric equation, but they are defined on the different parametric variable regions. In fact, $\mathcal{S}_i$ and $\mathcal{S}_r$ should be the different subregions contained in a larger surface. According to the assumption on the parametric equation $\mathbf{s}(\eta^\alpha)$, the coordinate curves of $\{\eta^{\alpha}\}$ constitute a net of curvature lines on $\mathcal{S}_i$. By virtue of the variable change $\eta^{1}(\theta^1,\theta^2)$ and $\eta^{2}(\theta^1,\theta^2)$, another parametric equation of surface $\mathcal{S}_i$ can be obtained as follow
\begin{equation}
\mathcal{S}_i:\ \mathbf{s}^*(\theta^\alpha) = \mathbf{s}(\eta^{1}(\theta^1,\theta^2),\eta^{2}(\theta^1,\theta^2)),\quad (\theta^\alpha) \in \Omega_r,
\label{Eq:ParaESi2}
\end{equation}
which is defined on the parametric variable region $\Omega_r$ in the $\theta^1\theta^2$-plane. We choose $\mathcal{S}_i$ as the shape of the base surface of the shell in an intermediate configuration $\mathcal{K}_i$. The position vector $\mathbf{X}^{*}$ of a material point in $\mathcal{K}_i$ is set to be (cf. Eq. \eqref{Eq:2})
\begin{equation}
    \begin{aligned}
        \mathbf{X}^{*} = &\mathbf{s}^*(\theta^\alpha) + Z\mathbf{n}^*(\theta^\alpha),\\
        = & \mathbf{s}(\eta^1(\theta^\alpha),\eta^2(\theta^\alpha)) + Z\mathbf{n}(\eta^1(\theta^\alpha),\eta^2(\theta^\alpha)), \quad (\theta^\alpha) \in \Omega_r,\ \ 0\leq Z\leq 2h.
    \end{aligned}
\label{Eq:posiV}
\end{equation}

Based on the above results, we can write out the growth tensor that produces the shape change of the base surface of the shell from $\mathcal{S}_r$ to $\mathcal{S}$. As shown in Fig. \ref{Fig:ChangeVariables}, the whole deformation process is divided into two steps. In the first step, we consider the deformation of the shell from the reference configuration $\mathcal{K}_r$ to the intermediate configuration $\mathcal{K}_i$ (i.e., the shape change of the base surface from $\mathcal{S}_r$ to $\mathcal{S}_i$ ). Based on \eqref{Eq:ParaESi2} and \eqref{Eq:posiV}, it is known that the corresponding deformation gradient tensor should be given by
\begin{equation}
    \begin{aligned}
        \mathbb{F}_0 =& (\nabla\mathbf{X^*})\mathbb{U}^{-1} + \frac{\partial\mathbf{X^*}}{\partial Z} \otimes \mathbf{n}\\
        =& \frac{\partial\eta^1}{\partial\theta^1} \mathbf{g}_{1}(\eta^\alpha)\otimes\mathbf{g}^1(\theta^\alpha) + \frac{1-\kappa_1 Z}{1-\kappa_2 Z}\frac{\partial\eta^1}{\partial\theta^2} \mathbf{g}_{1}(\eta^\alpha)\otimes\mathbf{g}^2(\theta^\alpha) \\
        & + \frac{1-\kappa_2 Z}{1-\kappa_1 Z}  \frac{\partial\eta^2}{\partial\theta^1} \mathbf{g}_2(\eta^\alpha) \otimes \mathbf{g}^1(\theta^\alpha) +  \frac{\partial\eta^2}{\partial\theta^2} \mathbf{g}_2(\eta^\alpha) \otimes \mathbf{g}^2(\theta^\alpha) + \mathbf{n}(\eta^\alpha) \otimes \mathbf{n}(\theta^\alpha),
    \end{aligned}
    \label{Eq:FG0}
\end{equation}
In Eq. \eqref{Eq:FG0}, the covariant base $\{\mathbf{g}_1,\mathbf{g}_2,\mathbf{n}\}$ is evaluated at the position $\mathbf{s}(\eta^\alpha)$ on $\mathcal{S}_i$ and the contravariant base $\{\mathbf{g}^1,\mathbf{g}^2,\mathbf{n}\}$ is evaluated at the position $\mathbf{s}(\theta^\alpha)$ on $\mathcal{S}_r$. In the second step, we consider the deformation of the shell from the intermediate configuration $\mathcal{K}_i$ to the current configuration $\mathcal{K}_t$ (i.e., the shape change of the base surface from $\mathcal{S}_i$ to $\mathcal{S}_t$). As shown in Fig. \ref{Fig:ChangeVariables}, $\mathcal{S}_i$ and $\mathcal{S}_t$ possess the same parametric variable region $\Omega_r^{*}$ in the $\eta^1\eta^2$-plane. Besides that, the coordinate curves of $\{\eta^\alpha\}$ constitute the net of curvature lines on these two surfaces. Thus, the formulas \eqref{Eq:GF1} and \eqref{Eq:GF3} obtained in section \ref{Sec:3.1} should be applicable in this case. The growth tensor that can induce the shape change from $\mathcal{S}_i$ to $\mathcal{S}_t$ is then given by
\begin{equation}
    \begin{aligned}
        \mathbb{G}_1 =& \left[1+Z \left(\kappa_1^*-\frac{L}{E}\right)\right]\sqrt{\frac{E}{E^*}}\mathbf{g}_1(\eta^\alpha)\otimes\mathbf{g}^1(\eta^\alpha) \\
        &+ \left[1+Z \left(\kappa_2^*-\frac{N}{G}\right)\right]\sqrt{\frac{G}{G^*}}\mathbf{g}_2(\eta^\alpha)\otimes\mathbf{g}^2(\eta^\alpha) + \mathbf{n}(\eta^\alpha)\otimes\mathbf{n}(\eta^\alpha).
    \end{aligned}
    \label{Eq:G1}
\end{equation}
In Eq. \eqref{Eq:G1}, $\{E,G,L,N\}$ are the fundamental quantities of surface $\mathcal{S}$ calculated with the parametric equation $\mathbf{x}^*(\eta^\alpha)$. $\{E^*,G^*\}$ and $\{\kappa_1^*,\kappa_2^*\}$ are the fundamental quantities and principal curvatures, respectively, of surface $\mathcal{S}_i$ calculated with the parametric equation $\mathbf{s}(\eta^\alpha)$. It can be directly verified that
\begin{equation}
    \mathbb{G}_1\mathbb{F}_0 = \mathbb{Q}\mathbb{G},
    \label{Eq:QG}
\end{equation}
where $\mathbb{Q}$ is the rotation tensor
\begin{equation}
\mathbb{Q}=\mathbf{g}_1(\eta^\alpha)\otimes\mathbf{g}^1(\theta^\alpha)+\mathbf{g}_2(\eta^\alpha)\otimes\mathbf{g}^2(\theta^\alpha) + \mathbf{n}(\eta^\alpha)\otimes\mathbf{n}(\theta^\alpha),
\end{equation}
and
\begin{equation}
    \begin{aligned}
        \mathbb{G} = & \left[1+Z \left(\kappa_1^*-\frac{L}{E}\right)\right]\sqrt{\frac{E}{E^*}}\frac{\partial\eta^1}{\partial\theta^1} \mathbf{g}_{1}(\theta^\alpha)\otimes\mathbf{g}^1(\theta^\alpha) \\
        & + \left[1+Z \left(\kappa_1^*-\frac{L}{E}\right)\right]\sqrt{\frac{E}{E^*}}\left(\frac{1-\kappa_1 Z}{1-\kappa_2 Z}\right)\frac{\partial\eta^1}{\partial\theta^2} \mathbf{g}_{1}(\theta^\alpha)\otimes\mathbf{g}^2(\theta^\alpha) \\
        & + \left[1+Z \left(\kappa_2^*-\frac{N}{G}\right)\right]\sqrt{\frac{G}{G^*}}\left(\frac{1-\kappa_2 Z}{1-\kappa_1 Z}\right)  \frac{\partial\eta^2}{\partial\theta^1} \mathbf{g}_2(\theta^\alpha) \otimes \mathbf{g}^1(\theta^\alpha) \\
        & +  \left[1+Z \left(\kappa_2^*-\frac{N}{G}\right)\right]\sqrt{\frac{G}{G^*}}\frac{\partial\eta^2}{\partial\theta^2} \mathbf{g}_2(\theta^\alpha) \otimes \mathbf{g}^2(\theta^\alpha) + \mathbf{n}(\theta^\alpha) \otimes \mathbf{n}(\theta^\alpha).
    \end{aligned}
    \label{Eq:GT2}
\end{equation}
Tensor $\mathbb{G}$ given in \eqref{Eq:GT2} is just the growth tensor that can result in the shape change of the base surface of the shell from $\mathcal{S}_r$ to $\mathcal{S}$ in the general case, which is consistent with the growth tensor obtained in \eqref{Eq:GT1} for the special case.


\subsection{A theoretical scheme for shape-programming}
\label{sec:3.3}

\begin{figure}[htbp]
    \centering
    \includegraphics[width=0.6\textwidth]{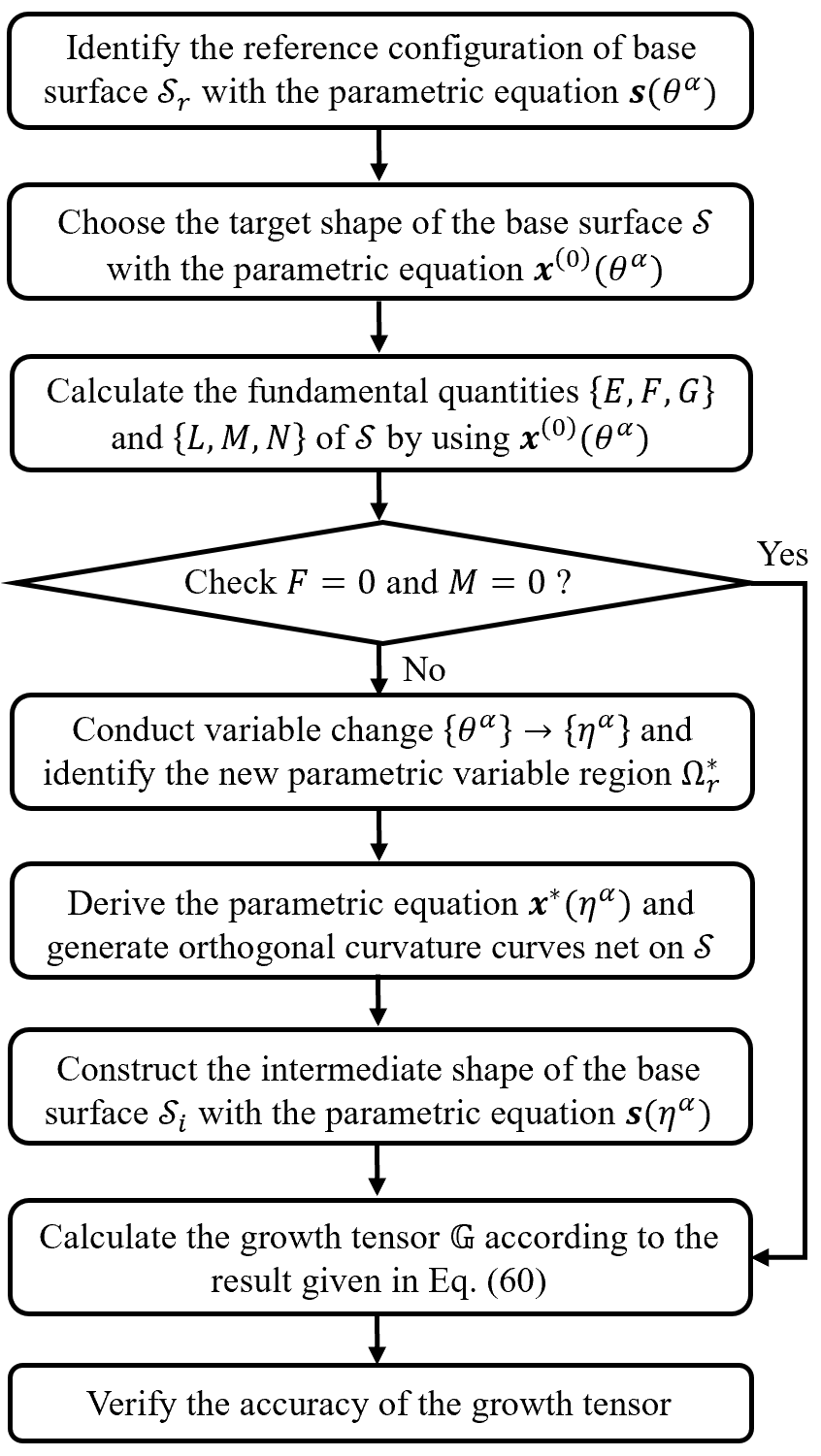}
    \caption{Flowchart of the theoretical scheme for shape-programming of a thin hyperelastic shell through differential growth.}
    \label{Fig:Scheme}
\end{figure}

Based on the above preparations, we propose a theoretical scheme for shape-programming of a thin hyperelastic shell through differential growth. The flowchart of this scheme is shown in Fig. \ref{Fig:Scheme}, which contains the following steps:

\begin{itemize}
\item With the given reference configuration $\mathcal{K}_r$ of the shell, we need to identify the parametric equation $\mathbf{s}(\theta^\alpha)$ for the initial shape of the base surface $\mathcal{S}_r$, which is defined on the region $\Omega_r$ of the $\theta^1\theta^2$-plane. By using $\mathbf{s}(\theta^\alpha)$, the fundamental quantities $\{E_r,G_r,L_r,N_r\}$ and the principal curvatures $\{\kappa_1,\kappa_2\}$ of surface $\mathcal{S}_r$ can be calculated.

\item We choose the target shape of the base surface $\mathcal{S}$, which has the parametric equation $\mathbf{x}^{(0)}(\theta^\alpha)$ defined on $\Omega_r$.

\item The fundamental quantities $\{E,F,G\}$ and $\{L,M,N\}$ of surface $\mathcal{S}$ are calculated by using the parametric equation $\mathbf{x}^{(0)}(\theta^\alpha)$. In the case $F=0$ and $M=0$, it is known that the parametric curves net of $\{\theta^\alpha\}$ is already an orthogonal net of curvature lines \citep{chenWH2017}. Then, the growth tensor $\mathbb{G}$ can be obtained from Eq. \eqref{Eq:GT1}.

\item If both $F$ and $M$ are not equal to zero, we need to conduct the variable change from $\{\theta^\alpha\}$ to $\{\eta^\alpha\}$, which yields a bijection from $\Omega_r$ to a new region $\Omega_r^*$ in the $\eta^1\eta^2$-plane. The explicit expressions of the variable change should be calculated from the differential forms given in \eqref{Eq:eta1eta2}, where the integrating factors $A_1^*$ and $A_2^*$ need to be determined in advance.

\item After the variable change, the surface $\mathcal{S}$ has a new parametric equation $\mathbf{x}^*(\eta^\alpha)$ defined on $\Omega_r^*$. The coordinate curves of $\{\eta^\alpha\}$ constitute an orthogonal net of curvature lines on $\mathcal{S}$.

\item By virtue of the variable change, an intermediate shape of the base surface $\mathcal{S}_i$ can be constructed, which has the parametric equation $\mathbf{s}(\theta^\alpha)$ defined on $\Omega_r$ and the parametric equation $\mathbf{s}^*(\eta^\alpha)$ defined on $\Omega_r^*$. The associated geometrical quantities of $\mathcal{S}_i$ are also calculated.

\item Based on the above results, the growth tensor $\mathbb{G}$ is calculated from Eq. \eqref{Eq:GT2}, which results in the shape change of the base surface of the shell from $\mathcal{S}_r$ to $\mathcal{S}$.

\item Finally, to check the correctness and accuracy of this scheme, the obtained growth tensor (or growth functions) is incorporated in a finite element analysis (we use Abaqus), and the growth-induced deformation of the shell is simulated.

\end{itemize}


\section{Examples}
\label{sec:4}

We demonstrate the feasibility and efficiency of the analytical framework for shape-programming of thin hyperelastic shells through differential growth proposed in Section \ref{sec:3} using some typical examples inspired by soft biological tissues in nature.

For the purpose of illustration, the reference configuration $\mathcal{K}_r$ of the shell is selected to be a cylindrical shell, which occupies the region $[R_0,R_0+2h]\times[0,\Theta_0]\times[0,l]$ within a cylindrical coordinate system in $\mathcal{R}^3$. The base face $\mathcal{S}_r$ of the shell has the following parametric equation
\begin{equation}
    \mathbf{s}(\theta^1,\theta^2) = \{ R_0 \cos(\theta^1), R_0 \sin(\theta^1), \theta^2 \}, \quad 0\leq\theta^1\leq\Theta_0,\ 0\leq\theta^2\leq l,
    \label{Eq:SrCylinder}
\end{equation}
where $\theta^1$ and $\theta^2$ are the parametric variables. It is clear that the coordinate curves of $\theta^1$ and $\theta^2$ constitute a net of curvature lines on $\mathcal{S}_r$. Besides that, we denote $Z = R - R_0$ ($R_0\leq R\leq R_0+2h$) as the thickness variable of the shell. From the parametric equation \eqref{Eq:SrCylinder}, we obtain the following covariant and contravariant base vectors on $\mathcal{S}_r$
\begin{equation}
    \begin{aligned}
    &\mathbf{g}_{1}=R_{0}\left[-\sin(\theta^1)\mathbf{e}_{1}+\cos(\theta^1)\mathbf{e}_{2}\right], \quad \mathbf{g}^{1} = \frac{\mathbf{g}_{1}}{R_{0}^{2}} , \\ &\mathbf{g}_{2} = \mathbf{g}^{2} = \mathbf{e}_{3}, \quad \mathbf{g}_{3}=\mathbf{g}^{3}=\mathbf{n}=\cos(\theta^1)\mathbf{e}_{1}+\sin(\theta^1)\mathbf{e}_{2}.
    \end{aligned}
    \label{basevec}
\end{equation}
The geometrical quantities of surface $\mathcal{S}_r$ are given by
\begin{equation}
    \begin{aligned}
        &E_r = R_0^2, \quad G_r = 1, \quad L_r = -R_0, \\
        &N_r = 0, \quad \kappa_1 = -1/R_0, \quad \kappa_2 = 0.
    \end{aligned}
    \label{geoq}
\end{equation}

\subsection{Example without change of variables}
In the first example, the target shape of the base surface $\mathcal{S}$ is selected to be a surface of revolution, which has the following parametric equation
\begin{equation}
    \mathbf{x}^{(0)}(\theta^1,\theta^2) = \{ u(\theta^2) \cos (\theta^1) , u(\theta^2) \sin (\theta^1) , v(\theta^2) \}, \quad 0\leq\theta^1\leq\Theta_0,\ 0\leq\theta^2\leq l,
    \label{Eq:SrCx}
\end{equation}
where $u(\theta^2)$ and $v(\theta^2)$ are arbitrarily smooth functions. Notice that both $\mathcal{S}_r$ and $\mathcal{S}$ have the parametric variable region $\Omega_r=[0,\Theta_0]\times[0,l]$. Corresponding to the parametric equation \eqref{Eq:SrCx}, the following first and second fundamental quantities of surface $\mathcal{S}$ are obtained
\begin{equation}
    \begin{aligned}
        & E = u^2, \quad F = 0, \quad G =u^{\prime 2} + v^{\prime 2}, \\
        & L = -\frac{u^2 v^{\prime} }{ \sqrt{ u^2 \left(u^{\prime 2} + v^{\prime 2}\right) }}, \quad M = 0, \quad N = \frac{u \left( v^{\prime} u^{\prime \prime} - u^{\prime} v^{\prime \prime}  \right) }{ \sqrt{ u^2 \left(u^{\prime 2} + v^{\prime 2}\right) }}.
    \end{aligned}
    \label{FQSR}
\end{equation}
Since $F=0$ and $M=0$, it is known that the $\theta^1$- and $\theta^2$-coordinate curves have already constituted an orthogonal net of curvature lines on $\mathcal{S}$. Therefore, the growth tensor in the shell should be set according to \eqref{Eq:GT1}, which contains the growth functions
\begin{equation}
    \begin{aligned}
        & \lambda_1^{(0)} = |u|, \quad  \lambda_2^{(0)} =  \sqrt{u^{\prime 2} + v^{\prime 2}} , \\
        & \lambda_1^{(1)} = -\frac{|u|}{R_0} + \frac{v^{\prime}}{  \sqrt{u^{\prime 2} + v^{\prime 2} }}, \quad
        \lambda_2^{(1)} = \frac{ u \left( u^{\prime} v^{\prime \prime} - v^{\prime} u^{\prime \prime}  \right) }{ |u| (u^{\prime 2} + v^{\prime 2})  }.
    \end{aligned}
    \label{Eq:GFRevo}
\end{equation}

For concreteness, we consider four kinds of revolution surfaces inspired by biological tissues, i.e., the sweet melon, the morning glory, the trachea and the apple. The parametric equations and the corresponding growth functions of these surfaces are listed in \eqref{Eq:ExampleSphere}-\eqref{Eq:ExampleApple}. To verify the accuracy of the obtained growth functions, we also conduct numerical simulations by using the UMAT subroutine in ABAQUS, where the constitutive relation of a compressible neo-Hookean material is adopted. The Poisson's ratio of the material is chosen to $\nu=0.4995$ to capture the effect of elastic incompressibility. The growth functions $\lambda_{1} = \lambda_{1}^{(0)}+Z\lambda_{1}^{(1)}$ and $\lambda_{2} = \lambda_{2}^{(0)}+Z\lambda_{2}^{(1)}$ are incorporated as the state variables in UMAT, which change gradually from $1$ to the target functions as those given in \eqref{Eq:ExampleSphere}-\eqref{Eq:ExampleApple}. The initial cylindrical shell has the dimensions $R_0=4$, $h=0.01$ and $l=4$. The value of $\Theta_0$ is set to $\pi$ or $2\pi$ depending upon the case. The whole sample is meshed into 20160 C3D8IH elements (8-node linear brick, hybrid, linear pressure, incompatible modes).

\begin{figure}[htbp]
    \centering
    \includegraphics[width=0.71\textwidth]{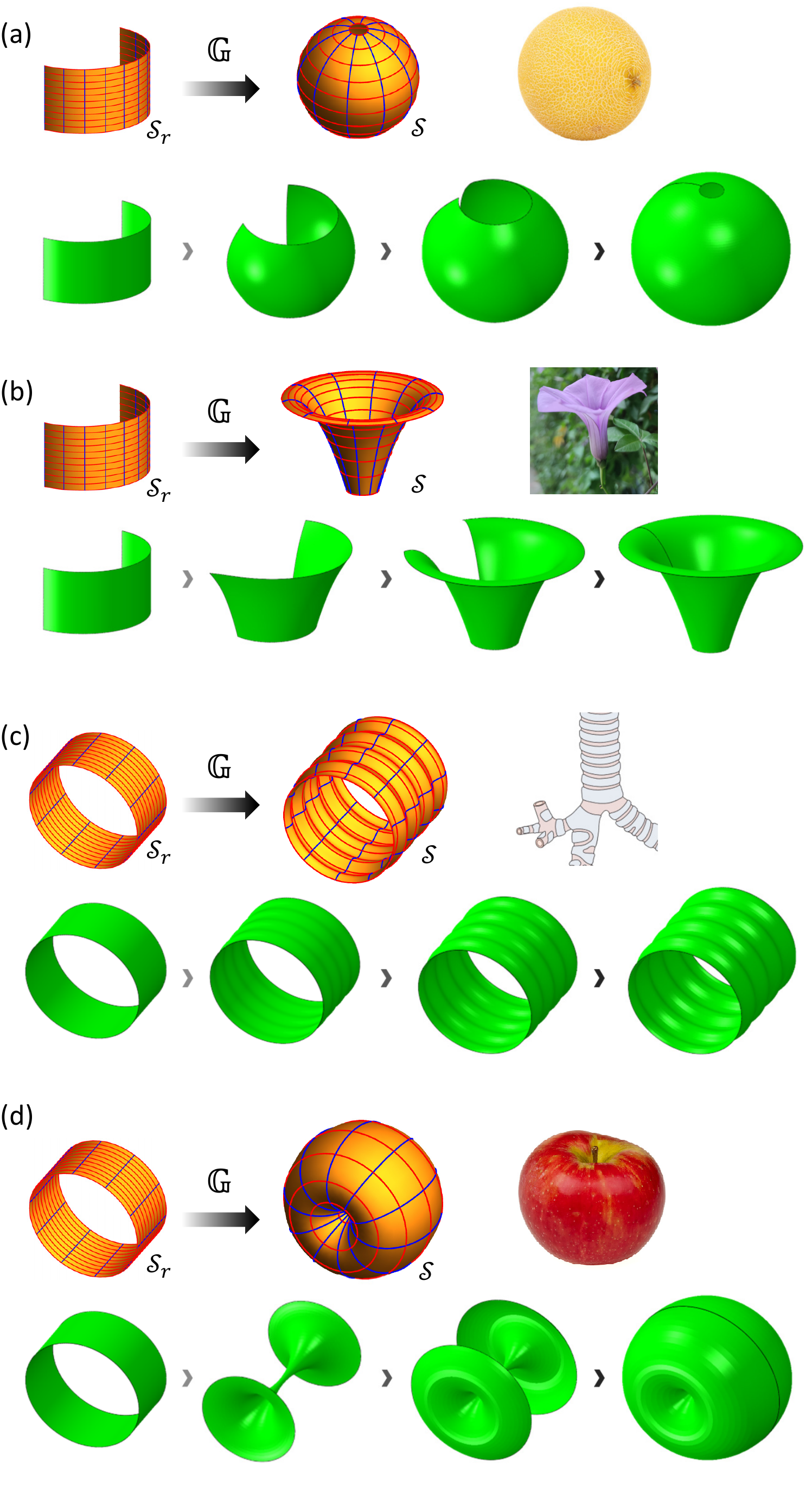}
    \caption{Numerical simulation results on the growing processes of the shells with the target surfaces and growth functions listed in \eqref{Eq:ExampleSphere}-\eqref{Eq:ExampleApple}: (a) the sweet melon;  (b) the morning glory; (c) the trachea; (d) the apple.}
    \label{Fig:Example1}
\end{figure}

In Fig. \ref{Fig:Example1}, we show the numerical simulation results. It can be seen that the grown states of the shells are in good agreement with the target shapes. Thus, the correctness of the obtained growth functions can be verified. It should be pointed out here that we only try to mimic the shapes of the different biological tissues, but we do not aim to reveal the underlying mechanisms responsible for the growth of the biological tissues.

\begin{itemize}
    \item Sweet melon $\left( 0 \leq \theta^1 \leq \pi, \quad 0 \leq \theta^2 \leq 4 \right):$
\end{itemize}
\begin{equation}
    \left\{
    \begin{aligned}
        & x^{(0)}= 4 \cos (2 {\theta^1}) \cos \left(\frac{9}{40} \pi  ({\theta^2}-2)\right), \\
        & y^{(0)}=4 \sin (2 {\theta^1}) \cos \left(\frac{9}{40} \pi  ({\theta^2}-2)\right), \\
        & z^{(0)}=-4 \cos \left(\frac{1}{40} \pi  (9 {\theta^2}+2)\right), \\
        & \lambda_{1} = 8 \cos \left(\frac{9}{40} \pi  ({\theta^2}-2)\right) , \quad \lambda_{2} =  \frac{9}{40} \pi  (Z+4) .
    \end{aligned}
    \right.
    \label{Eq:ExampleSphere}
\end{equation}
\begin{itemize}
    \item Morning glory $\left( 0 \leq \theta^1 \leq \pi, \quad 0 \leq {\theta^2} \leq 4 \right)$
\end{itemize}
\begin{equation}
    \left\{
    \begin{aligned}
        & x^{(0)}= -(1+{\theta^2}) \cos (2\theta^1), \\
        & y^{(0)}= -(1+{\theta^2}) \sin (2\theta^1), \\
        & z^{(0)}= 6-\frac{1}{8} \left( 7-2 {\theta^2} \right)^2, \\
        & \lambda_{1} = \frac{1}{2} ({\theta^2}+1) \left[4-\frac{4 (2 {\theta^2}-7) Z}{({\theta^2}+1)\sqrt{ 4 {\theta^2} \left( {\theta^2}-7\right)+53 }}-Z\right] ,\\
        &\lambda_{2} = \sqrt{{\theta^2}({\theta^2}-7) +\frac{53}{4}}-\frac{4 Z}{4 {\theta^2} \left( {\theta^2}-7\right)+53} .
    \end{aligned}
    \right.
\end{equation}
\begin{itemize}
    \item Trachea $\left( 0 \leq \theta^1 \leq 2\pi, \quad 0 \leq {\theta^2} \leq 4 \right)$
\end{itemize}
\begin{equation}
    \left\{
    \begin{aligned}
        & x^{(0)}= \frac{1}{5} \cos \theta^1  \left[20 + \sin(2 \pi {\theta^2}) \right], \\
        & y^{(0)}= \frac{1}{5} \sin \theta^1  \left[20 + \sin(2 \pi {\theta^2}) \right], \\
        & z^{(0)}= 2(2+{\theta^2}), \\
        & \lambda_{1} = -\frac{1}{20} (Z-4) \sin (2 \pi  {\theta^2}) + 4 + Z \left(\frac{5 \sqrt{2}}{\sqrt{\pi ^2 \cos (4 \pi  {\theta^2})+\pi ^2+50}}-1\right),\\
        & \lambda_{2} = \frac{1}{5} \sqrt{2} \sqrt{\pi^2 \cos (4 \pi  {\theta^2})+\pi^2+50}+\frac{20 \pi^2 Z \sin (2 \pi {\theta^2})}{\pi^2 \cos (4 \pi  {\theta^2})+\pi^2+50} .
    \end{aligned}
    \right.
\end{equation}
\begin{itemize}
    \item Apple $\left( 0 \leq \theta^1 \leq 2\pi, \quad 0 \leq {\theta^2} \leq 4 \right)$
\end{itemize}
\begin{equation}
    \left\{
    \begin{aligned}
        & x^{(0)}= 8 \cos \theta^1 \cos ^2\left(\frac{\pi {\theta^2}}{4}\right),  \\
        & y^{(0)}= 8 \sin \theta^1 \cos ^2\left(\frac{\pi {\theta^2}}{4}\right),  \\
        & z^{(0)}= -6 \sin \left(\frac{\pi  {\theta^2}}{2}\right),  \\
        & \lambda_{1} = \frac{-3 \sqrt{2} Z \cos \left(\frac{\pi  {\theta^2}}{2}\right)}{ \sqrt{5 \cos (\pi  {\theta^2})+13}} - 2(Z-4) \cos ^2\left(\frac{\pi  {\theta^2}}{4}\right) ,\\
        &\lambda_{2} = \pi  \left[\frac{\sqrt{5 \cos (\pi  {\theta^2})+13}}{\sqrt{2}}-\frac{6 Z}{5 \cos (\pi  {\theta^2})+13}\right] .
    \end{aligned}
    \right.
    \label{Eq:ExampleApple}
\end{equation}

\subsection{Example with change of variables}
To further demonstrate the efficiency of the proposed theoretical scheme, we study two more examples, in which the target shapes are chosen to be the Cereus Forbesii Spiralis and the tendril of pumpkin.

\begin{figure}[htbp]
    \centering
    \includegraphics[width=0.95\textwidth]{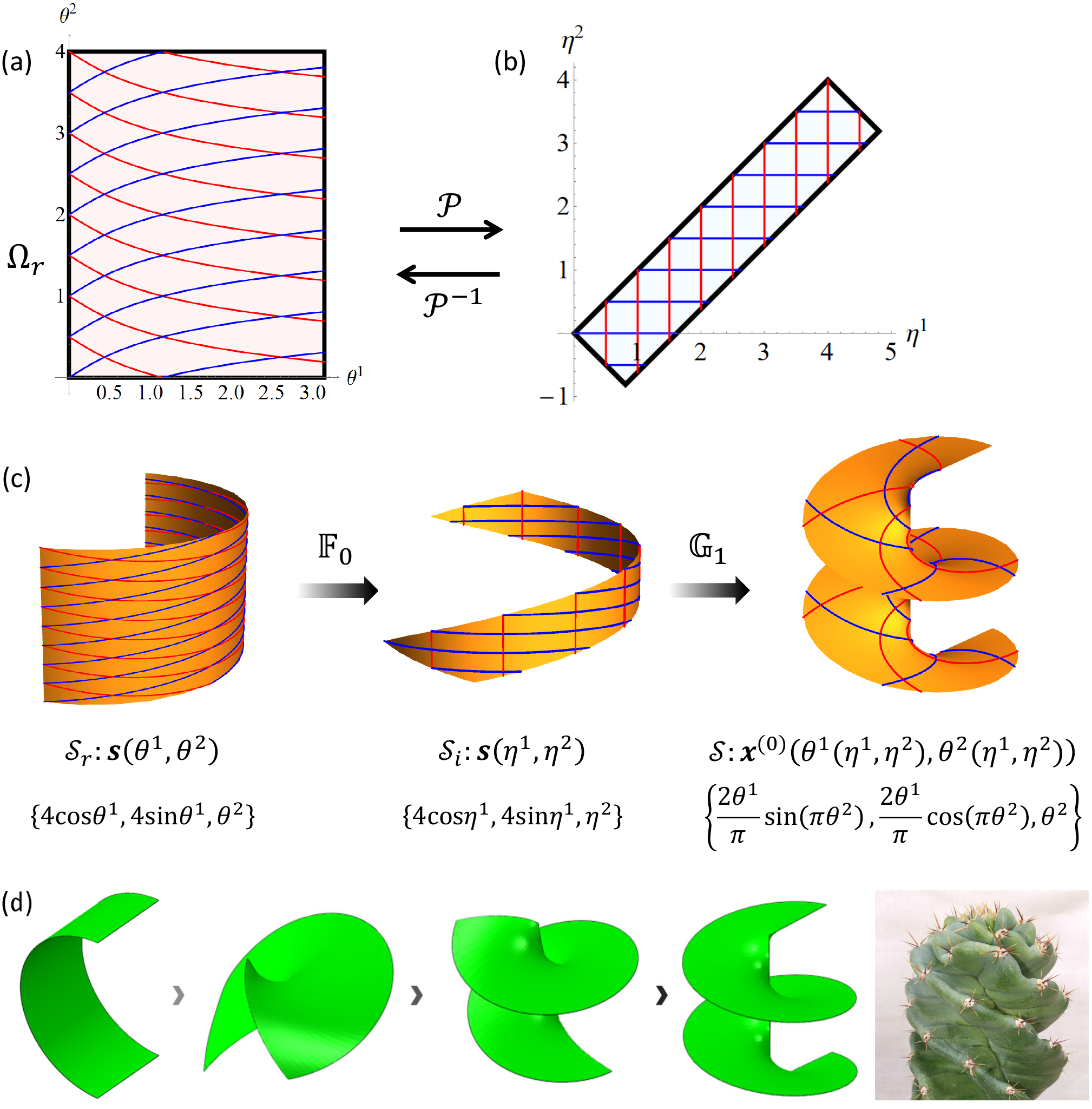}
    \caption{Change of variables between $\{ \theta^1,\theta^2 \}$ and $\{ \eta^1,\eta^2 \}$, and the decomposition of the growth process for generating Cereus Forbesii Spiralis: (a) the original region $\Omega_r$ on the parametric plane $\theta^1\theta^2$; (b) the new region $\Omega_r^*$ on the parametric plane $\eta^1\eta^2$; (c) the base surface $\mathcal{S}_r$ in the referential configuration $\mathcal{K}_r$, the base surface $\mathcal{S}_i$ in the intermediate configuration $\mathcal{K}_i$, and the target base surface $\mathcal{S}$ in the current configuration $\mathcal{K}_t$; (d) the simulated growing process for generating the Cereus Forbesii Spiralis configuration of the shell.}
    \label{Fig:MappingHelical}
\end{figure}

For the case of Cereus Forbesii Spiralis, the parametric equation of the target surface $\mathcal{S}$ is
\begin{equation}
    \mathbf{x}^{(0)}(\theta^1,\theta^2) = \left\{ \frac{2\theta^1}{\pi} \sin (\pi {\theta^2}) , \frac{2\theta^1}{\pi} \cos (\pi {\theta^2}), \theta^2 \right\},
    \label{Eq:ParaSpiral}
\end{equation}
where the region of the parametric variable is chosen to be $\Omega_r = [0, \pi]\times[0,4]$. Corresponding to this parametric equation, one can obtain the first and second fundamental quantities as follows
\begin{equation}
    \begin{aligned}
        & E = 4/\pi^2, \quad F = 0, \quad G = 1 + 4{\theta^1}^2, \\
        & L = 0, \quad M = 2/ \sqrt{ 1 + 4{\theta^1}^2}, \quad N = 0. \\
    \end{aligned}
    \label{Eq:EFGLMNSpiral}
\end{equation}
As the quantity $M \neq 0$, we need to conduct the change of variables from $\left(\theta^1,\theta^2 \right)$ to $\left(\eta^1,\eta^2 \right)$. According to the procedure of variable change introduced in section \ref{Sec:3.2}, we have
\begin{equation}
    \eta^1\left( \theta^1,{\theta^2} \right) = \frac{\text{arcsinh}(2\theta^1)}{\pi} + {\theta^2}, \quad \eta^2\left( \theta^1,{\theta^2} \right) = -\frac{\text{arcsinh}(2\theta^1)}{\pi} + {\theta^2} .
\end{equation}
After the variable transformation, the original region $\Omega_r$ in the $\theta^1\theta^2$- plane is mapped into a new region $\Omega_r^*$ in the $\eta^1\eta^2$-plane, which is shown in Fig. \ref{Fig:MappingHelical}(b).
On the region $\Omega_r^*$, a new surface $\mathcal{S}_i$ is defined as follows
\begin{equation}
    \mathcal{S}_i : \mathbf{s}^* = \{R_0 \cos \eta^1 , R_0 \sin \eta^1 , \eta^2 \}, \quad \eta^\alpha \in \Omega_r^*.
    \label{Eq:SiSpiral}
\end{equation}
Notice that the cylindrical shell $\mathcal{S}_r$ defined by \eqref{Eq:SrCylinder} and surface $\mathcal{S}_i$ defined by \eqref{Eq:SiSpiral} have the same parametric equation, but their parametric variable regions are different. As shown in Fig.\eqref{Fig:MappingHelical}c, both $\mathcal{S}_r$ and $\mathcal{S}_i$ can be viewed as a subregion cutting from a large cylindrical surface with radius $R_0=4$. Also, the coordinate curves of $\{\eta^\alpha\}$ (i.e., the blue and red curves in Fig.\ref{Fig:MappingHelical}) constitute the orthogonal nets of curvature lines on both $\mathcal{S}_i$ and $\mathcal{S}$. By choosing $\mathcal{S}_i$ as the base surface, we define an intermediate configuration $\mathcal{K}_i$ according to \eqref{Eq:posiV}, then the whole growth process can be divided into two steps:  $\mathcal{K}_r\rightarrow\mathcal{K}_i$ and from $\mathcal{K}_i\rightarrow\mathcal{K}_t$. For the first step, according to \eqref{Eq:FG0} the deformation gradient $\mathbb{F}_0$ is given by
\begin{equation}
    \begin{aligned}
        \mathbb{F}_0 = & \frac{2}{\pi  \sqrt{1+4{\theta^1}^2}} \mathbf{g}_{1}(\eta^\alpha)\otimes\mathbf{g}^1(\theta^\alpha) + (1+Z/R_0) \mathbf{g}_{1}(\eta^\alpha)\otimes\mathbf{g}^2(\theta^\alpha) \\
        & -\frac{2}{\pi \left(Z/R_0 + 1\right) \sqrt{1+4{\theta^1}^2}} \mathbf{g}_2(\eta^\alpha) \otimes \mathbf{g}^1(\theta^\alpha) +  \mathbf{g}_2(\eta^\alpha) \otimes \mathbf{g}^2(\theta^\alpha) + \mathbf{n}(\eta^\alpha) \otimes \mathbf{n}(\theta^\alpha).
    \end{aligned}
    \label{Eq:F0e1}
\end{equation}
For the second step, the growth tensor $\mathbb{G}_1$ on domain $\left(\eta^1,\eta^2 \right)$ is obtained according to \eqref{Eq:G1}
\begin{equation}
    \begin{aligned}
        \mathbb{G}_1 = & \mathbf{g}_{1}(\eta^\alpha)\otimes\mathbf{g}^1(\eta^\alpha) \Bigg[ \frac{1}{8} \sqrt{\cosh (\pi (\eta^1-\eta^2))+1} \\
        & -\frac{Z \left(16 \pi \sqrt{\cosh ^4\left(\frac{1}{2} \pi (\eta^1-\eta^2)\right)}+\cosh ^2(\pi (\eta^1-\eta^2))+2 \cosh (\pi (\eta^1-\eta^2))+1\right)}{32 (\cosh (\pi (\eta^1-\eta^2))+1)^{3/2}} \Bigg] \\
        &  +  \mathbf{g}_{2}(\eta^\alpha)\otimes\mathbf{g}^2(\eta^\alpha) \Big[ \frac{1}{2} \sqrt{\cosh (\pi (\eta^1-\eta^2))+1}    +\frac{2 \pi Z \sqrt{\cosh ^4\left(\frac{1}{2} \pi (\eta^1-\eta^2)\right)}}{(\cosh (\pi (\eta^1-\eta^2))+1)^{3/2}} \Big] \\
        & + \mathbf{n}(\eta^\alpha)\otimes\mathbf{n}(\eta^\alpha). \\
    \end{aligned}
    \label{Eq:G1e1}
\end{equation}
Then the tensor $\mathbb{G}$ generating the shape change from $\mathcal{S}_r$ to $\mathcal{S}$ can be obtained according to \eqref{Eq:GT2}.
To verify the correctness of these growth functions, we simulate the growth process of Cereus Forbesii Spiralis in ABAQUS.
The setting of the numerical simulations is the same as that introduced in the previous example.
The numerical results of this case are shown in Fig. \eqref{Fig:MappingHelical}d, which shows that the final shape of the shell can fit the target shape quite well.

\begin{figure}[htbp]
    \centering
    \includegraphics[width=0.95\textwidth]{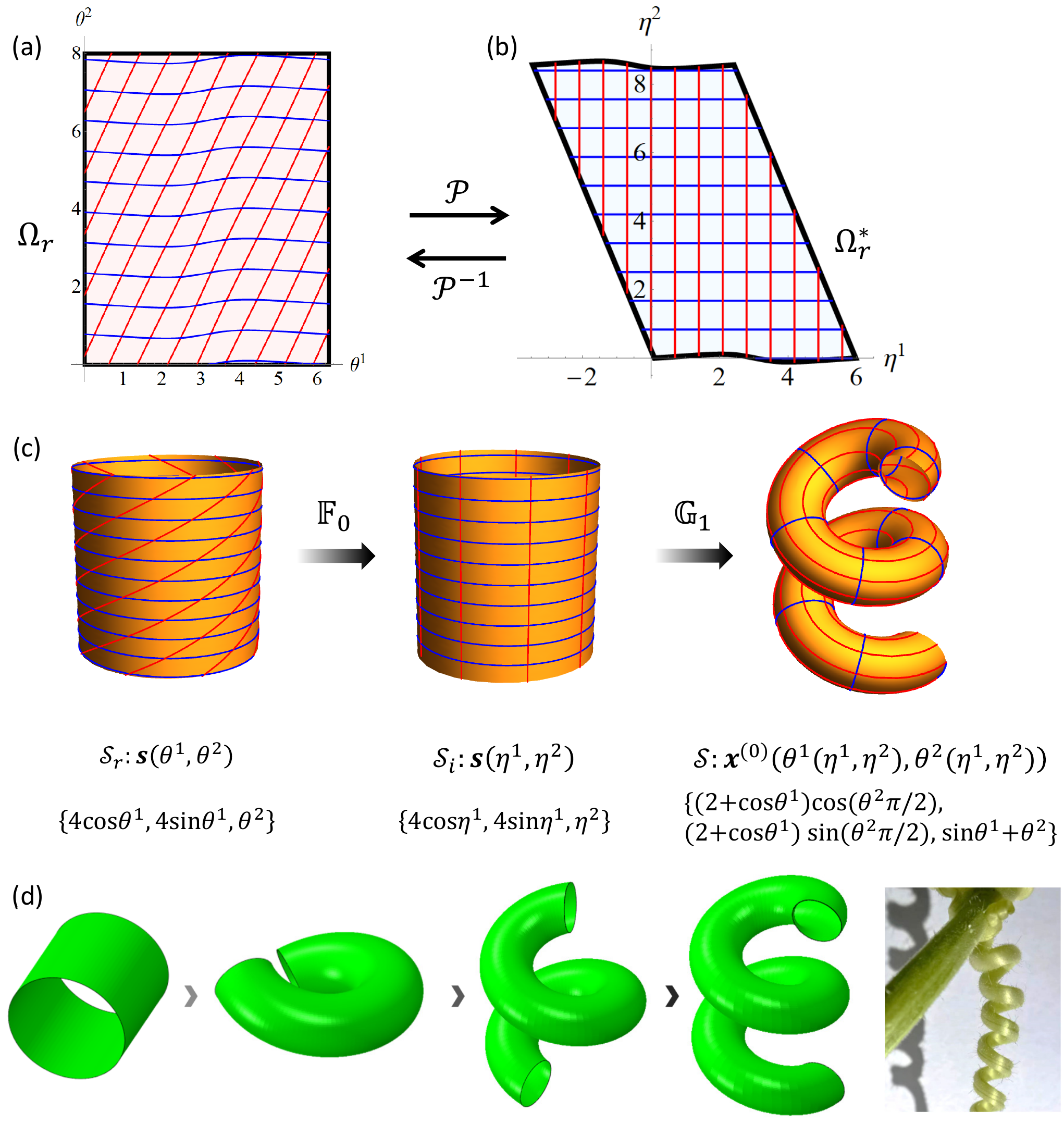}
    \caption{Change of variables between $\{ \theta^1,\theta^2 \}$ and $\{ \eta^1,\eta^2 \}$, and the decomposition of the growth process for generating tendril of pumpkin: (a) the original region $\Omega_r$ in the parametric plane $\theta^1\theta^2$; (b) the new region $\Omega_r^*$ in the parametric plane $\eta^1\eta^2$; (c) the base surface $\mathcal{S}_r$ in the referential configuration $\mathcal{K}_r$, the base surface $\mathcal{S}_i$ on the intermediate configuration $\mathcal{K}_i$, and the target base surface $\mathcal{S}$ in the current configuration $\mathcal{K}_t$; (d) the simulated growing process for generating the tendril of pumpkin configuration of the shell.}
    \label{Fig:MappingRod}
\end{figure}

For the case of pumpkin tendril, the parametric equation of the target surface $\mathcal{S}$ is
\begin{equation}
    \mathbf{x}^{(0)}(\theta^1,\theta^2) = \left\{ (\cos \theta^1 +2) \cos \left(\frac{\pi {\theta^2}}{2}\right), (\cos \theta^1 +2) \sin \left(\frac{\pi {\theta^2}}{2}\right), \sin \theta^1 +{\theta^2} \right\},
    \label{Eq:ParaTendril}
\end{equation}
where the region on the parametric plane is chosen to be $\Omega_r = [0, 2\pi] \times [0,8]$.
Corresponding to this parametric equation, one can obtain the first and second fundamental quantities as follow
\begin{equation}
    \begin{aligned}
        & E = 1, \quad F = \cos \theta^1, \quad G = \frac{1}{8} \left(\pi ^2 (8 \cos ({\theta^1})+\cos (2 {\theta^1}))+9 \pi ^2+8\right), \\
        & L = \frac{\sqrt{2} \pi  (\cos ({\theta^1})+2)}{\sqrt{8 \pi ^2 \cos ({\theta^1})+\left(\pi ^2-4\right) \cos (2 {\theta^1})+9 \pi ^2+4}}, \\
        & M = -\frac{\sqrt{2} \pi  \sin ^2({\theta^1})}{\sqrt{8 \pi ^2 \cos ({\theta^1})+\left(\pi ^2-4\right) \cos (2 {\theta^1})+9 \pi ^2+4}}, \\
        & N = \frac{\pi ^3 \cos ({\theta^1}) (\cos ({\theta^1})+2)^2}{2 \sqrt{2} \sqrt{8 \pi ^2 \cos ({\theta^1})+\left(\pi ^2-4\right) \cos (2 {\theta^1})+9 \pi ^2+4}}. \\
    \end{aligned}
    \label{Eq:EFGLMNRod}
\end{equation}
Note that the quantities $F \neq 0$ and $M \neq 0$, thus the change of variables from $\left(\theta^1,\theta^2 \right)$ to $\left(\eta^1,\eta^2 \right)$ is required.
As shown in Fig.\ref{Fig:MappingRod}(b), the original region $\Omega_r$ is mapped into a new region $\Omega_r^*$ through the change of variables.
Following the same parametric equation, these two regions define surfaces $\mathcal{S}_r$ and $\mathcal{S}_i$ respectively, where $\mathcal{S}_r$ is a cylinder with radius $R_0 = 4$ and length $l = 8$, while $\mathcal{S}_i$ is an irregular shaped subregion cut from a cylinder with radius $R_0 = 4$ as follow
\begin{equation}
    \mathcal{S}_i : \mathbf{s}^* = \{R_0 \cos \eta^1 , R_0 \sin \eta^1 , \eta^2 \}, \quad \eta^\alpha \in \Omega_r^*.
    \label{Eq:SiRod}
\end{equation}
Also, the coordinate curves $\{\eta^\alpha\}$(i.e., the blue and red curves in Fig.\ref{Fig:MappingRod}) constitute an orthogonal curvature net on both $\mathcal{S}_i$ and $\mathcal{S}$.
By choosing $\mathcal{S}_i$ as the base surface, we define an intermediate configuration $\mathcal{K}_i$.
The whole shape morphing process can be divided into two steps: from $\mathcal{K}_r$ to $\mathcal{K}_i$ described by $\mathbb{F}_0$, and the from $\mathcal{K}_i$ to $\mathcal{K}_t$ induced by $\mathbb{G}_1$.
However, the analytical explicit expressions for integrating factors and $\left(\eta^1,\eta^2 \right)$ are difficult to obtain.
Therefore the tensor $\mathbb{F}_0$ and $\mathbb{G}_1$ are calculated numerically in this case, and then the growth values are passed to the relating integration points on meshes in ABAQUS.
According to the numerical results of this case shown in Fig.\eqref{Fig:MappingRod}(d), the final shapes fit the target shapes quite well.


\section{Conclusions}

The large deformations of thin hyperelastic shells induced by differential growth were investigated in this paper. To fulfill the goal of shape-programming of thin hyperelastic shells, the following tasks have been accomplished: (i) a consistent finite-strain shell equation system for modeling the growth-induced deformations of incompressible hyperelastic shells was formulated; (ii) the problem of shape-programming was solved analytically under the stress-free condition, from which the explicit expressions of the growth functions in terms of the geometrical quantities (i.e., the first and second fundamental forms) of the target surfaces were derived; (iii) a general theoretical scheme for shape-programming of thin hyperelastic shells through differential growth was proposed; (iv) to verify the correctness and efficiency of the scheme, some typical examples were studied, where the configurations of some biological tissues were simulated with the obtained growth functions.

Since the formulas derived in the current work have relatively simple forms and are valid for general incompressible hyperelastic material, the presented theoretical scheme for shape-programming would have wide potential applications for design and manufacturing of intelligent soft devices. Furthermore, the analytical results can also shed light on understanding the mechanical behaviors of some soft biological tissues in nature during the growth process. It should be noted that the current work still has some shortcomings that need to be tackled in future. One shortcoming is that the explicit expressions of growth tensor is derived based on the stress-free assumption, therefore, they are not valid for shell samples subjected to external loads or boundary restrictions. Additionally, the theoretical scheme is not applicable to complex 3D surfaces without explicit parametric equations. In that case, an efficient numerical scheme for shape-programming of complicated surfaces needs to be proposed.

\section*{Supplementary material}
Movie 1: Growth process of shape-programming cases. Video of the growing processes of the six illustrative examples introduced in Fig. \ref{Fig:Example1}, \ref{Fig:MappingHelical} and \ref{Fig:MappingRod} of the main text, which is available at \url{https://github.com/Jeff97/growth-deformation-of-shell}


\section*{Declaration of competing interest}

The authors declare that they have no known competing financial interests or personal relationships that could have appeared to influence the work reported in this paper.


\section*{Acknowledgments}
\label{sec:ack}

This work is supported by the National Natural Science Foundation of China (Project No.: 11872184). Z.L. is supported by the China Scholarship Council (CSC) Grant \#202106150121. M.H. and Z. L. are indebted to the funding through an Engineering and Physical Sciences Research Council (EPSRC) Impact Acceleration Award (EP/R511614/1).


\appendix

\section{Some results for incompressible Neo-Hookean material}
\label{app:1}
\setcounter{equation}{0}
\renewcommand\theequation{A.\arabic{equation}}

To obtain some concrete results on the unknowns ($\mathbf{x}^{(n)}$ and $p^{(n)}$), we further assume that the shell is made of neo-Hookean material with the following elastic strain-energy function
\begin{equation}
    \phi(\mathbb{F},\mathbb{G})=J_{G}\phi_0(\mathbb{A})=J_GC_0\left[\mathrm{tr}(\mathbb{A}\mathbb{A}^T)-3\right],
\end{equation}
where $C_{0}$ is a material constant. From the elastic strain-energy function $\phi(\mathbb{F},\mathbb{G})$, the nominal stress tensor $\mathbb{S}$ is given by
\begin{equation}
    \mathbb{S}=J_{G}\mathbb{G}^{-1}\left(2C_{0}\mathbb{A}^{T}-p(R,Z)\mathbb{A}^{-1}\right).
    \label{Eq:Constitutive}
\end{equation}

For simplicity, we assume the shell is under traction-free condition.
By taking series expansion on $Z=0$ and through some truncation manipulation, a closed linear system for $\{ \mathbf{x}^{(1)}, \mathbf{x}^{(2)}, p^{(0)}, p^{(1)} \}$ is formulated by \eqref{Eq:ExpanDetA} and \eqref{Eq:ExpanDivS}$_1$, combining with the boundary conditions \eqref{Eq:Bou3}$_1$.
Then the following expressions of $\{ \mathbf{x}^{(1)}, \mathbf{x}^{(2)}, p^{(0)}, p^{(1)} \}$ in terms of $\mathbf{x}^{(0)}$ are solved
\begin{equation}
    \mathbf{x}^{(1)} = \Lambda \frac{\mathbf{x}_{N}}{\Delta}, \quad p^{(0)}=2 C_0 \frac{\Lambda^2}{\Delta},
\end{equation}
\begin{equation}
    \begin{aligned}
         \mathbf{x}^{(2)} =& \frac{1}{{{\Delta^{5/2}}}}{\mathbf{x}_N}\left( {{\Lambda^2}{t_9} + \frac{{{\Delta^2}}}{{{\Lambda^2}}}{t_8} - {\Delta^{3/2}}{t_1}} \right) - \frac{1}{{{\Lambda^2}}}{\mathbf{a}} \\
         &+ \frac{1}{{{\Delta^3}{\Lambda^3}}}\left[ {{\mathbf{x}^{(0)}_{,1}}\left( {{\Lambda^4}{t_5} - {\Delta^3}\lambda_{2}^{(0)} {t_7}} \right) + {\mathbf{x}^{(0)}_{,2}}\left( {{\Lambda^4}{t_4} - {\Delta^3}\lambda_{1}^{(0)} {t_6}} \right)} \right], \\
        p^{(1)}  =& 2C_0 \left( {\frac{1}{{\Lambda \sqrt \Delta }}{t_8} - \frac{{2\Lambda }}{\Delta}{t_1} + \frac{{{\Lambda^3}}}{{{\Delta^{5/2}}}}{t_9}} \right),
    \end{aligned}
    \label{Eq:r2p1-NeoHookean}
\end{equation}

where
\begin{equation*}
    \begin{aligned}
      & \Lambda = \lambda_{1}^{(0)} \lambda_{2}^{(0)} ,\quad \mathbf{x}_{N} = \mathbf{x}^{(0)}_{,1} \times \mathbf{x}^{(0)}_{,2} ,\quad \Delta = \mathbf{x}_{N} \cdot \mathbf{x}_{N}, \\
      & B_{\alpha \beta}=\mathbf{g}_\alpha \cdot \mathbf{g}_{\beta,\alpha}, \quad \mathbf{a} = \left( {\lambda_{1}^{(0)}}^2 \mathbf{x}^{(0)}_{,2,2} + {\lambda_{2}^{(0)}}^2 \mathbf{x}^{(0)}_{,1,1} \right), \\
      & {t_1} = \left( {{\kappa_{1}} + {\kappa_{2}}} \right)\Lambda - \lambda_{1}^{(1)}\lambda_{2}^{(0)} - \lambda_{1}^{(0)}\lambda_{2}^{(1)}, \\
      & {t_2} = \left( {{B_{11}} + {B_{21}}} \right)\Lambda ,\quad {t_3} = \left( {{B_{12}} + {B_{22}}} \right)\Lambda , \\
      & {t_4} = \Lambda \left( {{\Delta_{,1}}F - {\Delta_{,2}}E} \right) + \Delta \left[ {E\left( {2{\Lambda_{,2}} + {t_3}} \right) - F\left( {2{\Lambda_{,1}} + {t_2}} \right)} \right], \\
      & {t_5} = \Lambda \left( {{\Delta_{,2}}F - {\Delta_{,1}}G} \right) + \Delta \left[ {G\left( {2{\Lambda_{,1}} + {t_2}} \right) - F\left( {2{\Lambda_{,2}} + {t_3}} \right)} \right], \\
      & {t_6} = \lambda_{1,2}^{(0)} \Lambda  + \lambda_{1}^{(0)}\left( {{t_3} - \lambda_{2,2}^{(0)}\lambda_{1}^{(0)}} \right), \quad {t_7} = \lambda_{2,1}^{(0)} \Lambda  + \lambda_{2}^{(0)}\left( {{t_2} - \lambda_{1,1}^{(0)} \lambda_{2}^{(0)}} \right), \\
      & {t_8} = \left( \lambda_{1}^{{(0)}^2} N  + \lambda_{2}^{{(0)}^2} L  \right), \quad {t_9} = \left( {E N  - 2F M  + G L } \right).
    \end{aligned}
\end{equation*}

The expressions of $\mathbb{S}^{(0)}$ in terms of $\mathbf{x}^{(0)}$ are obtained by substituting $\mathbf{x}^{(1)}$ and $p^{(0)}$ into \eqref{Eq:S0S1JMPS}$_1$
\begin{equation}
    \begin{aligned}
        \mathbb{S}^{(0)} =& 2C_0 \mathbf{g}_1 \otimes \left[ \frac{\Lambda^3}{\Delta^2} (F \mathbf{x}^{(0)}_{,2} - G \mathbf{x}^{(0)}_{,1}) + \frac{\lambda_{2}^{(0)}}{\lambda_{1}^{(0)}} \mathbf{x}^{(0)}_{,1} \right]\\
        &+ 2C_0\mathbf{g}_2 \otimes \left[ -\frac{\Lambda^3}{\Delta^2} (F \mathbf{x}^{(0)}_{,1} - E \mathbf{x}^{(0)}_{,2}) + \frac{\lambda_{1}^{(0)}}{\lambda_{2}^{(0)}} \mathbf{x}^{(0)}_{,2} \right].\\
    \end{aligned}
    \label{Eq:S0-Neo}
\end{equation}
Accordingly, expression of $\mathbb{S}^{(1)}$ is also obtained, where $\mathbf{x}^{(2)}$ and $p^{(1)}$ are kept for brevity
\begin{equation}
    \begin{aligned}
        \mathbb{S}^{(1)} =&2C_0 \mathbf{g}_1 \otimes \Bigg[ \frac{{\left( { - \lambda_1^{(1)}\lambda_2^{(0)} + \lambda_1^{(0)}\lambda_2^{(1)} + \Lambda {\kappa_{1}}} \right)}}{\lambda_{1}^{(0)^2}}{\mathbf{x}^{(0)}_{,1}} - \frac{{{\Lambda^2}}}{\Delta}{\mathbf{x}^{(0)}_{, 2}} \times {\mathbf{x}^{(2)}}\\
        &+ \frac{{{\Lambda^4}}}{{{\Delta^{5/2}}}}\left( {N{\mathbf{x}^{(0)}_{,1}} - M{\mathbf{x}^{(0)}_{, 2}}} \right) + \frac{{\lambda_2^{(0)}\left( {\Delta {\Lambda_{,1}} - \Lambda  {\Delta_{,1}}} \right)}}{{{\Delta^2} \lambda_1^{(0)}}}{\mathbf{x}_N}\\
        &- \frac{{\Lambda \left( { {\kappa_{2}}{\Lambda^2} + {p^{(1)}}\Delta/(2{C_0})} \right)}}{{{\Delta^2}}}\left( {G{\mathbf{x}^{(0)}_{,1}} - F{\mathbf{x}^{(0)}_{, 2}}} \right) + \frac{ \lambda_{2}^{(0)^2}}{\Delta} \mathbf{x}_{N,1} \Bigg]\\
        &+ 2C_0\mathbf{g}_2 \otimes \Bigg[ \frac{{\left( {\lambda _1^{(1)}\lambda _2^{(0)} - \lambda _1^{(0)}\lambda _2^{(1)} + \Lambda {\kappa_{2}}} \right)}}{{\lambda _2^{{{(0)}^2}}}}{\mathbf{x}^{(0)}_{,2}} + \frac{{{\Lambda ^2}}}{\Delta }{\mathbf{x}^{(0)}_{,1}} \times  {\mathbf{x}^{(2)}}\\
        &- \frac{{{\Lambda ^4}}}{{{Q^{5/2}}}}\left( {M{\mathbf{x}^{(0)}_{,1}} - L{\mathbf{x}^{(0)}_{,2}}} \right) + \frac{{\lambda _1^{(0)}\left( {\Delta {\Lambda _{,2}} - \Lambda {\Delta _{,2}}} \right)}}{{{\Delta ^2}\lambda _2^{(0)}}}{\mathbf{x}_N}\\
        &+ \frac{{\Lambda \left( {{\kappa_{1}}{\Lambda^2} + {p^{(1)}}\Delta /(2{C_0})} \right)}}{{{\Delta ^2}}}\left( {F{\mathbf{x}^{(0)}_{,1}} - E{\mathbf{x}^{(0)}_{,2}}} \right) + \frac{{\lambda _1^{{{(0)}^2}}}}{\Delta }{\mathbf{x}_{N,2}} \Bigg]\\
        &+ 2C_0 \mathbf{n} \otimes \Bigg[ \left( { - \frac{{\Lambda {t_1}}}{\Delta } - {p^{(1)}}/(2{C_0})} \right){\mathbf{x}_N} + \Lambda \mathbf{x}^{(2)}\\
        &+ \frac{{{\Lambda^2}}}{{{\Delta^3}}}\left[ {\left( {\Delta {\Lambda _{,1}} - \Lambda {\Delta _{,1}}} \right)\left( {G{\mathbf{x}^{(0)}_{,1}} - F{\mathbf{x}^{(0)}_{,2}}} \right) - \left( {\Delta {\Lambda _{,2}} - \Lambda {\Delta _{,2}}} \right)\left( {F{\mathbf{x}^{(0)}_{,1}} - E{\mathbf{x}^{(0)}_{,2}}} \right)} \right]\\
        &+ \frac{{{\Lambda^3}}}{{{\Delta^2}}}\left( {{\mathbf{x}_{N,2}} \times {\mathbf{x}^{(0)}_{,1}} - {\mathbf{x}_{N,1}} \times {\mathbf{x}^{(0)}_{,2}}} \right) \Bigg].\\
    \end{aligned}
    \label{Eq:S1-Neo}
\end{equation}
Note that $\mathbb{S}^{(1)}$ is also in terms of $\mathbf{x}^{(0)}$ with the use of \eqref{Eq:r2p1-NeoHookean}.

In order to fulfil the goal of shape-programming, growth functions $\{ \lambda_{1}^{(0)},\lambda_{1}^{(1)},\lambda_{2}^{(0)},\lambda_{2}^{(1)} \}$ of an arbitrary target shape $\mathbf{x}^{(0)}$ need to be determined from shell equation system.
Generally, the growth functions of a certain target shape may not be unique.
To facilitate derivation, we assume all the components  $\mathbb{S}^{(0)}$ and $\mathbb{S}^{(1)}$ in current configuration $\mathcal{K}_t$ are zero.
It is clear that, under the stress-free assumption, shell equation \eqref{Eq:ShellEqs} and boundary conditions \eqref{Eq:ShellBcs} are satisfied automatically.

First, all components of $\mathbb{S}^{(0)}$ in \eqref{Eq:S0-Neo} are set to be zero
\begin{equation}
\left\{
    \begin{aligned}
        & \frac{\Lambda^3}{\Delta^2} (F \mathbf{x}^{(0)}_{,2} - G \mathbf{x}^{(0)}_{,1}) + \frac{\lambda_{2}^{(0)}}{\lambda_{1}^{(0)}} \mathbf{x}^{(0)}_{,1} = \mathbf{0}, \\
        - & \frac{\Lambda^3}{\Delta^2} (F \mathbf{x}^{(0)}_{,1} - E \mathbf{x}^{(0)}_{,2}) + \frac{\lambda_{1}^{(0)}}{\lambda_{2}^{(0)}} \mathbf{x}^{(0)}_{,2} = \mathbf{0}. \\
    \end{aligned}
\right.
    \label{Eq:S0-Neo2}
\end{equation}
For simplicity, we assume $F=0$ in the current configuration, which means the moving frame $\{ \mathbf{x}^{(0)}_{,1},\mathbf{x}^{(0)}_{,2},\mathbf{x}_{N} \}$ are perpendicular to each other. Then the equations \eqref{Eq:S0-Neo2} are simplified as
\begin{equation}
\left\{
    \begin{aligned}
        &   \left(-G \frac{\Lambda^3}{E^2 G^2} + \frac{\lambda_{2}^{(0)}}{\lambda_{1}^{(0)}}\right)\mathbf{x}^{(0)}_{,1} = \mathbf{0}, \\
        &  \left(-E \frac{\Lambda^3}{E^2 G^2} + \frac{\lambda_{1}^{(0)}}{\lambda_{2}^{(0)}}\right) \mathbf{x}^{(0)}_{,2} = \mathbf{0}, \\
    \end{aligned}
\right.
\end{equation}
where the relation $\Delta=E G$ is used.
Subsequently, growth functions $\lambda_{1}^{(0)}$ and $\lambda_{2}^{(0)}$ are solved
\begin{equation}
    \lambda_{1}^{(0)}=\sqrt{E}, \quad \lambda_{2}^{(0)}=\sqrt{G},
    \label{Eq:GF1Neo}
\end{equation}
where growth functions $\lambda_{1}^{(0)}$ and $\lambda_{2}^{(0)}$ just represent extension or shrinkage along the coordinate curves $\{ \theta^\alpha \}$ on $\mathcal{S}_r$.

Second, we consider all components of $\mathbb{S}^{(1)}$ in \eqref{Eq:S1-Neo} are zero
\begin{equation}
\left\{
    \begin{aligned}
        &\frac{{\left( { - \lambda_1^{(1)}\lambda_2^{(0)} + \lambda_1^{(0)}\lambda_2^{(1)} + \Lambda {\kappa_{1}}} \right)}}{\lambda_{1}^{(0)^2}}{\mathbf{x}^{(0)}_{,1}} - \frac{{{\Lambda^2}}}{\Delta}{\mathbf{x}^{(0)}_{, 2}} \times {\mathbf{x}^{(2)}}\\
        & + \frac{{{\Lambda^4}}}{{{\Delta^{5/2}}}}\left( {N{\mathbf{x}^{(0)}_{,1}} - M{\mathbf{x}^{(0)}_{, 2}}} \right) + \frac{{\lambda_2^{(0)}\left( {\Delta {\Lambda_{,1}} - \Lambda  {\Delta_{,1}}} \right)}}{{{\Delta^2} \lambda_1^{(0)}}}{\mathbf{x}_N}\\
        & - \frac{{\Lambda \left( { {\kappa_{2}}{\Lambda^2} + {p^{(1)}}\Delta/(2{C_0})} \right)}}{{{\Delta^2}}} G{\mathbf{x}^{(0)}_{,1}}   + \frac{ \lambda_{2}^{(0)^2}}{\Delta} \mathbf{x}_{N,1} =\mathbf{0}, \\
        &  \frac{{\left( {\lambda _1^{(1)}\lambda _2^{(0)} - \lambda _1^{(0)}\lambda _2^{(1)} + \Lambda {\kappa_{2}}} \right)}}{{\lambda _2^{{{(0)}^2}}}}{\mathbf{x}^{(0)}_{,2}} + \frac{{{\Lambda ^2}}}{\Delta }{\mathbf{x}^{(0)}_{,1}} \times  {\mathbf{x}^{(2)}}\\
        & - \frac{{{\Lambda ^4}}}{{{Q^{5/2}}}}\left( {M{\mathbf{x}^{(0)}_{,1}} - L{\mathbf{x}^{(0)}_{,2}}} \right) + \frac{{\lambda _1^{(0)}\left( {\Delta {\Lambda _{,2}} - \Lambda {\Delta _{,2}}} \right)}}{{{\Delta ^2}\lambda _2^{(0)}}}{\mathbf{x}_N}\\
        & - \frac{{\Lambda \left( {{\kappa_{1}}{\Lambda^2} + {p^{(1)}}\Delta /(2{C_0})} \right)}}{{{\Delta ^2}}}  E{\mathbf{x}^{(0)}_{,2}} + \frac{{\lambda _1^{{{(0)}^2}}}}{\Delta }{\mathbf{x}_{N,2}} = \mathbf{0},\\
        & \left( { - \frac{{\Lambda {t_1}}}{\Delta } - {p^{(1)}}/(2{C_0})} \right){\mathbf{x}_N} + \Lambda \mathbf{x}^{(2)}\\
        & + \frac{{{\Lambda^2}}}{{{\Delta^3}}}\left[ {\left( {\Delta {\Lambda _{,1}} - \Lambda {\Delta _{,1}}} \right)G{\mathbf{x}^{(0)}_{,1}}  + \left( {\Delta {\Lambda _{,2}} - \Lambda {\Delta _{,2}}} \right) E{\mathbf{x}^{(0)}_{,2}}} \right]\\
        & + \frac{{{\Lambda^3}}}{{{\Delta^2}}}\left( {{\mathbf{x}_{N,2}} \times {\mathbf{x}^{(0)}_{,1}} - {\mathbf{x}_{N,1}} \times {\mathbf{x}^{(0)}_{,2}}} \right) = \mathbf{0}.\\
    \end{aligned}
\right.
    \label{Eq:S1Free}
\end{equation}
With the use of \eqref{Eq:GF1Neo} and $\Delta = EG$, \eqref{Eq:S1Free}$_3$ is automatically satisfied and \eqref{Eq:S1Free}$_1$ and \eqref{Eq:S1Free}$_2$ have the following form
\begin{equation}
\left\{
    \begin{aligned}
        & - \left( {EN + G\left( {2L - E\left( {2 \kappa_1 +  \kappa_2} \right) + 2\sqrt E \lambda_{1}^{(1)}} \right) + E\sqrt G \lambda_{2}^{(1)}} \right)  \mathbf{x}^{(0)}_{,1}  = M E {\mathbf{x}^{(0)}_{,2}} ,\\
        & - \left( {2EN + G\left( {L - E\left( { \kappa_1 + 2 \kappa_2} \right) + \sqrt E \lambda_{1}^{(1)}} \right) + 2E\sqrt G \lambda_{2}^{(1)}} \right) { \mathbf{x}^{(0)}_{,2} } = M E {\mathbf{x}^{(0)}_{,1}} .\\
    \end{aligned}
\right.
    \label{Eq:S1Free1}
\end{equation}
To ensure the holds of Eqs. \eqref{Eq:S1Free1}, we need to set $M = 0$, which together with $F=0$ assume that the coordinate curves $\{ \theta^\alpha \}$ formulate the orthogonal net of curvature lines on the target surface $\mathcal{S}_t$. Subsequently, the growth functions $\lambda_{1}^{(1)}$ and $\lambda_{2}^{(1)}$ are solved
\begin{equation}
    \lambda_1^{(1)} = \left(\kappa_1 - \frac{L}{E}\right)\sqrt{E},\quad \lambda_2^{(1)} = \left(\kappa_2 - \frac{N}{G}\right)\sqrt{G}.
    \label{Eq:GF2Neo}
\end{equation}
It can be seen that the growth functions \eqref{Eq:GF1Neo} and \eqref{Eq:GF2Neo} are coincident with the results obtained in Section \ref{Sec:3.1}.
Compared with the plate sample in \citet{Wang2022}, a distinct feature of the current growth functions is that, the effects of curvature $\kappa_1$ and $\kappa_2$ are taken into account. By solving the problem of shape-programming of the Neo-Hookean shell, the relations between growth functions and geometric properties of the base surface $\mathcal{S}$ are also revealed.

\bibliographystyle{unsrtnat} 
\bibliography{refs}

\end{document}